\documentstyle[12pt]{article}
\newtheorem{theorem}{Theorem}[section]
\newtheorem{prop}{Proposition}[section]

\newcommand{\U}{U_q(X_N^{(r)})}

\renewcommand{\a}{\alpha}

\newcommand{\la}{\lambda}
\newcommand{\be}{\beta}
\newcommand{\ka}{\kappa}
\newcommand{\ep}{\epsilon}

\title{Vertex Operators for Twisted Quantum Affine Algebras}

\author{Naihuan Jing
and Kailash C. Misra}
\begin{document}           

\maketitle                 
\begin{center}
Department of Mathematics\\
North Carolina State University\\
Raleigh, NC 27695-8205\\
USA
\end{center}
\date{June 28, 1996}

\begin{abstract}
We construct explicitly the $q$-vertex operators (intertwining operators)
for the level one
modules $V(\Lambda_i)$ of the classical quantum affine algebras
of twisted types using interacting bosons, where $i=0, 1$ for $A_{2n-1}^{(2)}$,
$i=0$ for $D_4^{(3)}$,
$i=0, n$ for $D_{n+1}^{(2)}$, and $i=n$ for $A_{2n}^{(2)}$.
A perfect crystal graph for $D_4^{(3)}$ is constructed as a by-product.
\end{abstract}

\section{Introduction}
  In a major step toward understanding $q$-conformal field theory
\cite{kn:FR}, Frenkel
and Reshetikhin introduced the $q$-Knizhnik-Zamoldchikov equation
associated with the quantum affine algebras. The main theoretical tool they
used is that of certain intertwining operators called $q$-vertex operators
between two types of representations of quantum affine algebras. They showed
that the matrix coefficients of these vertex operators will give the
solutions of the $q$-KZ equations.

  On the other hand, in series of works by Kyoto School and their
collaborators elsewhere (cf. the monograph
\cite{kn:JM} for details) we see that the very same vertex operators
take the role of corner transfer matrices in vertex models in
statistical mechanics. The subsequent
analysis enabled them to provide rigorous mathematical theory to explain the
appearance of Baxter's method of the corner transfer matrix. The space of
physical states are
understood as a tensor product of highest weight representations of
quantum affine algebras \cite{kn:DFJMN}. In level one case since the explicit
realization of
the underlying quantum
affine algebra \cite{kn:FJ} is available, the method works particularly well.
For instance, one feature of their work is that the correlation
functions can be computed in an explicit form using the bosonic
realizations of the $q$-vertex operators \cite{kn:JMMN}.

  The program of understanding various vertex models relies upon the explicit
knowledge of the $q$-vertex operators associated quantum affine algebras.
So far the explicit realization of level one $q$-vertex operators are carried
out for the quantum affine algebras in the cases of $A_n^{(1)}$ \cite{
kn:JMMN, kn:K},
$D_n^{(1)}$ \cite{kn:JKK} and $B_n^{(1)}$ \cite{kn:JM}.
The level two cases of $A_1^{(1)}$ were treated in \cite{kn:I} and,
 as special cases of our general construction, in \cite{kn:JM}
. Some other cases for reducible
modules are also given in \cite{kn:KQS, kn:M}.

  In this paper we will give explicit construction of level one
  $q$-vertex operators for twisted
quantum affine algebras of types $A_{2n+1}^{(2)}, D_{n+1}^{(2)}$ ,
$A_{2n}^{(2)}$ and $D_4^{(3)}$. Our work is based on the
explicit realization of level one
twisted quantum affine algebras given first in \cite{kn:J1} for $\Lambda_0$.
We generalize the realization to other level one irreducible modules in
Section 3.

 In section two we will review some preliminary results about quantum affine
algebras and Drinfeld realization, which is written in a modified
form \cite{kn:J3}. We then recall the level one realization of
twisted quantum affine algebras \cite{kn:J1} and in particular we extend the
construction to all other level one modules by the coset method. In the final
section we give the explicit construction of the level one $q$-vertex operators
for classical twisted quantum affine algebras.

In the previous construction of $q$-vertex operators, the finite dimensional
level zero representations are exactly coming from the crystal graphs of
the quantum affine algebras. We found that this is no longer true for
the twisted cases. In the case of $D_4^{(3)}$ we constructed a
$8$-dimensional representation and its associated perfect crystal
graph (c.f. (\ref{E:D_4})). As far as we know, this is the first
perfect crystal structure
for $D_4^{(3)}$.

Another novelty is that the $q$-Heisenburg subalgebra in
the Drinfeld generators of the level zero modules are no longer
diagonalizable. This is due to the underlying vector space representation
may admit a special $1$-dimensional subspace, which are used in the
perfect crystal graphs for types $D_{n+1}^{(2)}$, $A_{2n}^{(2)}$ and
$D_4^{(3)}$ (c.f. Section 2).

N. J. acknowledges partial support by NSA grants MDA 904-94-H-2061 and
MDA 904-96-1-0087. K. M. acknowledges partial support by NSA grant
MDA 904-96-1-0013.

\section{Twisted Quantum affine algebras $\U$}

Let $\bf g$ be the finite dimensional simple Lie algebra ${\bf{sl}}(2n)$
($n\geq 3$)
, ${\bf{so}}(2n+2)$ ($n\geq 2$), ${\bf{sl}}(2n+1)$
($n\geq 1$) or $\bf{so}(8)$. We associate an integer $r$ to the four
cases as follows: $r=2, 2, 2, 3$ respectively. We denote the Dynkin diagram
by $\Gamma$ and label their simple
roots and coroots as follows:
$$\a_1', \cdots, \a_N'; \hskip 20pt h_1', \cdots, h_N'$$
where $N=2n-1, n+1, 2n, 4$ respectively, $\a_i\in {\bf h}'^*$
and ${\bf h}'={\bf C}h_1'\oplus \cdots\oplus {\bf C}h_N'$
is the Cartan subalgebra of $\bf g$. We denote their
Chevalley generators by
$$e_i', f_i', h_i', \ i=1, \cdots, N.$$
The root and weight
lattices are then defined by:
\begin{eqnarray*}
Q'&=&{\bf Z}\a_1'+\cdots+ {\bf Z}\a_N' \\
P'&=&{\bf Z}\la_1'+\cdots +{\bf Z}\la_N'
\end{eqnarray*}
where $\la'_i$ are the fundamental weights of the Lie algebra $\bf g$:
$\la_i'(h_j')=\delta_{ij}, i,j= 1, \cdots N$.

We let $<.|.>$ be the nondegenerate invariant bilinear form on $\bf g$
normalized by  $<\a|\a>=2$, i.e.,
$<x|y>=tr(xy), \frac12tr(xy), tr(xy), \frac 12tr(xy)$
respectively. Then $(h_i'|h_i')=2$ for $i=1, \cdots, N$. Since
Lie algebra
$\bf g$ is simply-laced, we can identify the invariant form on $\bf h'^*$
to that of $\bf h'$.

Let $\sigma$ be a diagram automorphism of order $r=2, 3$ for the Lie algebra
$\bf g$ such that
\begin{eqnarray}
&&\sigma(h_i')=h_{N-i}', i=1, \cdots N, \ \ \mbox{for type } \ A_{2n-1}
\mbox{ or } \ A_{2n},\nonumber \\
&&\sigma(h_i')=h_i', i=1, \cdots , n-1=N-2; \sigma(h_n')=h_{n+1}',
\ \ \mbox{for type} \ D_{n+1} \nonumber\\
&&\sigma(h'_1, h'_2, h'_3, h'_4)=(h'_3, h_2', h'_4, h'_1)\ \ \mbox{for type} \
D_4^{(3)}\nonumber
\end{eqnarray}
Then the Lie algebra is decomposed as a ${\bf Z}/r{\bf Z}$-graded
Lie algebra:
$${\bf g}={\bf g}_0\oplus \cdots\oplus{\bf g}_{r-1},$$
where ${\bf g}_i=\{ x\in {\bf g}| \sigma(x)=\omega^i x\}$ and $
\omega=e^{2\pi\sqrt{-1}/r}$.
The Dynkin diagram $\Gamma$ decomposes
itself into $\sigma$-orbits: $\Gamma=\Gamma_0 \cup \cdots \cup\Gamma_{r-1}$.
It is well-known that
the subalgebra ${\bf g}_0$ is the simple Lie algebra of types
$C_n$, $B_n$, $B_n$ and $G_2$ respectively.

The twisted affine Lie algebras are realized as a subalgebra of the central
extension of the loop algebra:
$$\hat{{\bf g}}(\sigma)=\sum_i {\bf g}_{i \ mod \ r}\otimes t^i \oplus {\bf C}c
\oplus {\bf C}d,$$
where $c$ is the central element, $d$ the degree element, and the
Lie bracket is given by
$$
[x\otimes t^i, y\otimes t^j]=[xy]\otimes t^{i+j}+\frac ir\delta_{i, -j}
<x|y>c.
$$
The twisted affine Lie algebra $\hat{{\bf g}}(\sigma)$ is also
generated by the Chevalley generators:
\begin{eqnarray*}
e_i=E_i\otimes t^{\delta_{i0}}, \ f_i=F_i\otimes
t^{-\delta_{i0}}, h_i=H_i\otimes 1+ra_{\ep}^{-1}c,
\ \mbox{for } i=0, 1, \cdots, n.
\end{eqnarray*}
where $a_{\ep}=a_0=1$ unless ${\bf g}={\bf sl}(2n+1)$ in which
$a_{\ep}=a_n=2$. The elements $E_i, F_i, H_i$ ($i=0, 1\cdots n$)
are defined as follows:
\begin{eqnarray*}
&&E_i=e_i', F_i=f_i, h_i=h_i', \ \mbox{if } \sigma(i)=i; \\
&&E_i=\frac 1r\sum_{j=0}^{r-1}e'_{\sigma^j(i)}, \
F_i=\sum_{j=0}^{r-1}f'_{\sigma^j(i)}, \
H_i=\frac 1r\sum_{j=0}^{r-1}h'_{\sigma^j(i)}, \ \mbox{if } \sigma(i)\neq i, \\
&&E_n=\frac 1{2}(e'_n+e'_{n+1}),
F_n=(f'_n+f'_{n+1}),
H_n=(h'_n+h'_{n+1}) \ \mbox{for } A_{2n}, r=2;\\
&&E_0=\frac 1r\sum_{j=0}^{r-1}\omega^jf'_{-\sigma^j(\theta^0)},
F_0=-\sum_{j=0}^{r-1}\omega^je_{\sigma^j(\theta^0)},
H_0=-\frac 1r\sum_{j=0}^{r-1}h'_{-\sigma^j(\theta^0)}, \  \mbox{except for }
A_{2n}, r=2\\
&&E_0=\frac 1rf'_{-\theta^0}, F_0=-e'_{\theta^0}, H_0=-\frac 1rh'_{\theta^0}, \
\mbox{for }
A_{2n}, r=2;\\
&&\theta^0=\left\{\begin{array}{ll}
\a'_1+\cdots+\a'_{2n-2}, \ \mbox{for } A_{2n-1}, r=2\\
\a'_1+\cdots+\a'_n , \ \mbox{for } D_{n+1}, r=2\\
\a'_1+\cdots+\a'_{2n}, \ \mbox{for } A_{2n}, r=2 \\
\a'_1+\a'_2+\a'_3, \ \mbox{for } D_4, r=3
\end{array}\right.
\end{eqnarray*}

The above Chevalley generators give a realization of the twisted affine
Kac-Moody
algebra $\hat{\bf g}(A)$ associated to the affine Cartan matrices \cite{kn:K}
$A=(A_{ij}), i,j\in I=\{0, 1, \cdots, n\}$:
\begin{eqnarray}
&&A_{2n-1}^{(2)}
=2\sum_{i=0}^nE_{ii}-\sum_{i=1}^{n-1}(E_{i,i+1}+E_{i+1,i})-E_{02}-E_{20}-
E_{n-1,n},\nonumber\\
&&D_{n+1}^{(2)}=
2\sum_{i=0}^nE_{ii}-\sum_{i=0}^{n-1}(E_{i,i+1}+E_{i+1,i})-E_{01}-E_{n,n-1}
\nonumber\\
&&A^{(2)}_{2n}
=2\sum_{i=0}^nE_{ii}-\sum_{i=0}^{n-1}(E_{i,i+1}+E_{i+1,i})-E_{10}-
(1+\delta_{n,1})E_{n,n-1}.
\nonumber\\
&&D_4^{(3)}=2\sum_{i=0}^2E_{ii}-\sum_{i=0}^2(E_{i,i+1}+E_{i+1,i})-
2E_{12}\nonumber
\end{eqnarray}
where $E_{ij}$'s are the unit matrices in ${\bf Z}^{(n+1)\times (n+1)}$.
We remark that our choice of $A^{(2)}_{2n}$ differs from Kac's convention in
\cite{kn:K}. Moreover we define that
$$
\a_i=\frac 1p_i\sum_{s=0}^{r-1}\a_{\sigma^s(i)}',
$$
where $p_i=1$ if $\sigma(i)\neq i$ and $p_i=r$ if $\sigma(i)=i$, and extended
to
the weight lattice. It is
necessary to normalize the invariant form on ${\bf g}_0$ by
\begin{equation}\label{E:form}
(\a, \beta)=<\a|\beta>/r.
\end{equation}

We define the affine root lattice of $\hat{\bf g}$ to be
$$\hat{Q}={\bf Z}\alpha_0\oplus\cdots\oplus{\bf Z}\alpha_n,
$$
where $\alpha_i\in \bf h^*$, and $<\alpha_i, h_j>=A_{ji}$, $<\alpha_i, d>=0$.
The affine weight lattice $\hat{P}$ is defined to be
$$
\hat{P}={\bf Z} \Lambda_0\oplus\cdots\oplus{\bf Z} \Lambda_n\oplus{\bf Z}
\delta ,
$$
where
$
\Lambda_i(h_j)=\delta_{ij}, \Lambda_i(d)=0$, and
$\delta(h_j)=0, \delta(d)=1 $ for $j\in I$.
The dual affine weight lattice is then defined as
$$
\hat{P}^{\vee}={\bf Z} h_0 \oplus {\bf Z} h_1
\oplus \cdots \oplus {\bf Z} h_n \oplus {\bf Z} d .
$$

The symmetric bilinear form $<\ |\ >/r$ on ${\bf h}'^*$ induces
the normalized symmetric form on ${\bf h}^*$ denoted as $(\ |\ )$ (c.f.
(\ref{E:form}))
satisfying that
\begin{equation}
(\alpha_i|\alpha_j)=d_ia_{ij}, \ \ (\delta|\alpha_i)=(\delta|\delta)=0
\ \ \mbox{for all} \ i,j\in I \label{E:2.2},
\end{equation}
where $(d_0, \cdots, d_n)=(1, \cdots, 1, 2)$, $(1, 2, \cdots, 2, 1)$,
$(2, 1, \cdots, 1, 1/2)$, and $(1, 1, 3)$,
for $A_{2n-1}^{(2)}$, $D_{n+1}^{(2)}$, $A_{2n}^{(2)}$ and $D_4^{(3)}$
respectively.

Let $q_i=q^{d_i}=q^{\frac 12(\a_i|\a_i)}, i\in I$.
The quantum affine algebra $U_q(X_N^{(r)})$ is
 the associative algebra with 1 over ${\bf C}(q^{1/2})$
generated by the elements $e_i$, $f_i$ $(i\in I)$ and $q^h$ $(h\in
\hat{P}^{\vee})$
with the following relations :
\begin{eqnarray}
\ q^h q^{h'}&=&q^{h+h'} \ \ \mbox{for} \ h,h'\in \hat{P}^{\vee},
\quad q^0=1,
\nonumber\\
q^h e_{i} q^{-h}&=&q^{\alpha_i(h)} e_{i}, \ \
q^h f_{i} q^{-h}=q^{-\alpha_i(h)} f_{i} \ \ \mbox{for} \ h\in \hat{P}^{\vee}
(i\in I), \nonumber\\
e_{i}f_{j}-f_{j}e_{i} &=&\delta_{ij}
\displaystyle\frac {t_i-t_i^{-1}} {q_i-q_i^{-1}}, \ \mbox{where} \
t_i=q_i^{h_i}=q^{\frac 12(\a_i|\a_i)h_i} \  \mbox{and} \ i,j \in I,
\label{E:2.5}\\
\sum_{m+k=1-a_{ij}}&& (-1)^m e_{i}^{(m)} e_{j} e_{i}^{(n)}=0, \quad
\mbox{and}\nonumber\\
\sum_{m+n=1-a_{ij}}&& (-1)^m f_{i}^{(m)} f_{j} f_{i}^{(n)}=0
 \ \ \mbox{for} \  i\neq j,  \nonumber
\end{eqnarray}
where $e_{i}^{(k)}=e_{i}^k/[k]_i !$, $f_{i}^{(k)}=f_{i}^k/[k]_i !$,
$[m]_i!=\prod_{k=1}^m [k]_i$, and
$[k]_i=\displaystyle\frac {q_i^k-q_i^{-k}} {q_i-q_i^{-1}}.$
For simplicity we will write $[k]_i=[k]$ for $q_i=q$.
The derived subalgebra generated by $e_i$, $f_i$, $t_i$ $(i \in I)$ is
denoted by $U'_q(X_N^{(r)})$.

The algebra $\U$ has a
Hopf algebra structure with comultiplication
$\Delta$, counit $\varepsilon$, and antipode $S$ defined by
\begin{eqnarray}  \nonumber
\Delta(q^h)&=&q^h \otimes q^h \ \ \mbox{for} \ h\in \hat{P}^{\vee},\\
\Delta(e_{i})&=&e_{i}\otimes 1 + t_i \otimes e_{i},
\nonumber\\
\Delta(f_{i})&=&f_{i}\otimes t_i^{-1} + 1 \otimes f_{i} \ \ \mbox{for}
\ i\in I,  \nonumber \\
\varepsilon (q^h)&=&1 \ \ \mbox{for} \ h\in \hat{P}^{\vee}
,\label{E:2.6}\\
\varepsilon(e_{i})&=&\varepsilon(f_{i})=0
\ \ \mbox{for} \ i\in I, \nonumber\\
S(q^h)&=&q^{-h} \ \ \mbox{for} \ h\in \hat{P}^{\vee}, \nonumber\\
S(e_{i})&=&-t_i^{-1}e_{i}, \ \ S (f_{i})=-f_{i}t_i \ \
\mbox{for} \ i\in I.   \nonumber
\end{eqnarray}

Let $V, W$ be two $\U$-modules. The tensor product $V\otimes W$ is defined
as the $\U$-module via the coproduct $\Delta$.
The (restricted) dual $\U$-module $V^*$ is defined by
$$(x \cdot v^*) (u)=v^*(S(x)\cdot u) $$
for $x\in U_q({\bf g})$, $u\in V$, and $v^*\in V^*$.

We now recall Drinfeld's realization of the quantum affine
algebra $\U$ (and of $U_q'(X_N^{(r)})$) \cite{kn:D, kn:J3}. We will
present a slightly different form \cite{kn:J3}
to avoid the $h$-adic completion.

Let $\omega$ be a primitive $r$th root. We denote $[k]_j=
\frac{q^k_i-q^{-k}_i}{q_i-q^{-1}_i}$ if $j$ belongs to the $\sigma$-orbit of
$i$, then $[k]_j$ is defined for all $j=1, \cdots, N$ though we frequently
use only $[k]_i$ for $i\in \{0\}\cup\Gamma_{\sigma}=\{0, 1, \cdots, n\} $.

\begin{theorem} \cite{kn:D, kn:J3}
Let ${\bf U}$ be the associative algebra with 1 over ${\bf C}(q^{1/2})$
generated by the elements $x_i^{\pm}(k)$, $a_i(l)$, $K_i^{\pm 1}$,
$\gamma^{\pm 1/2}$, $q^{\pm d}$ $(i=1,2,\cdots,N, k\in {\bf Z},
l\in {\bf Z} \setminus \{0\})$ with the following defining relations :
\begin{eqnarray}
&&x^{\pm}_{\sigma(i)}(k)=\omega^k x^{\pm}_i(k),  \ \
a_{\sigma(i)}(l)=\omega^l a_i(l) ,\nonumber\\
&&\mbox{} [\gamma^{\pm 1/2}, u]=0 \ \ \mbox{for all} \ u\in \textstyle {\bf U},
\nonumber\\
&&\mbox{} [a_i(k), a_j(l)]=\delta_{k+l,0}
\displaystyle\sum_{s=0}^{r-1}\frac {[k<\a'_i|\sigma^s(\a'_j)>/d_i]_i}{k}
\omega^{ks}\displaystyle\frac {\gamma^k-\gamma^{-k}}{q_j-q^{-1}_j},
\nonumber\\
&&\mbox{} [a_i(k), K_j^{\pm 1}]=[q^{\pm d}, K_j^{\pm 1}]=0,
\nonumber\\
&&q^d x_i^{\pm}(k) q^{-d}=q^k x_i^{\pm } (k), \ \
q^d a_i(l) q^{-d}=q^l a_i(l),
\nonumber\\
&&K_i x_j^{\pm}(k) K_i^{-1}=q^{\pm (\alpha_i|\alpha_j)} x_j ^{\pm}(k),
\nonumber\\
&&\mbox{} [a_i(k), x_j^{\pm} (l)]=\pm \displaystyle\sum_{s=0}^{r-1}
\frac {[k<\a'_i|\sigma^s(\a'_j)>/d_i]_i}{k}\omega^{ks}
\gamma^{\mp |k|/2} x_j^{\pm}(k+l),
\nonumber\\
&&\prod_s(z-\omega^sq^{\pm<\a'_i|\sigma^s(\a'_j)>)}w)x_i^{\pm}(z)x^{\pm}_j(w)
=\prod_s(zq^{\pm<\a'_i|\sigma^s(\a'_j)>}-\omega^s w)
x_j^{\pm}(w)x^{\pm}_i(z)
\nonumber\\
&&\mbox{} [x_i^{+}(k), x_j^{-}(l)]=\displaystyle
\sum_{s=0}^{r-1}\frac {\delta_{\sigma^s(i), j}\omega^{sl}}{q_i-q_i^{-1}}
\left( \gamma^{\frac {k-l}{2}} \psi_{i} (k+l)-\gamma^{\frac{l-k}{2}}
\varphi_{i} (k+l) \right),
\label{E:2.11}
\end{eqnarray}
where $\psi_{i}(m)$ and $\varphi_{i}(-m)$ $(m\in {\bf Z}_{\ge 0})$
are defined by
\begin{eqnarray}
&&\sum_{m=0}^{\infty} \psi_{i}(m) z^{-m}
=K_i \textstyle {exp} \left( (q_i-q_i^{-1}) \sum_{k=1}^{\infty} a_i(k) z^{-k}
\right),
\nonumber\\
&&\sum_{m=0}^{\infty} \varphi_{i}(-m) z^{m}
=K_i^{-1} \textstyle {exp} \left(- (q_i-q_i^{-1}) \sum_{k=1}^{\infty} a_i(-k)
z^{k}
\right),
\nonumber\\
\mbox{Sym}_{z_1, z_2}&&P_{ij}^{\pm}(z_1, z_2)\sum_{s=0}^{2}(-1)^s
\left[\begin{array}{c} 2\\s\end{array}\right]_{q^{d_{ij}}}
x^{\pm}_{i}(z_1)\cdots
x^{\pm}_i(z_s)x^{\pm}_j(w)x^{\pm}_i(z_{s+1})\cdots x^{\pm}_i(z_2)=0,
\nonumber\\
&&\hskip 1in \mbox{for } A_{ij}=-1, \sigma(i)\neq j,\nonumber\\
\ \mbox{Sym}_{z_1, z_2, z_3}&&\left[(q^{\mp 3r/4}z_1-{q^{r/4}+q^{-r/4})z_2+
q^{\pm 3r/4}z_3)x_i^{\pm}(z_1)}x_i^{\pm}(z_2)x_i^{\pm}(z_3)\right]=0,
\nonumber\\
&&\hskip 1in \mbox{for } A_{i,\sigma(i)}=-1\nonumber\\
\end{eqnarray}
where Sym means the symmetrization over
  $ z_i $,
  $ P_{ij}^{\pm}(z,w) $
and
  $ d_{ij} $
are defined as follows:

\begin{eqnarray}
&&\mbox{If }
  \sigma(i) = i,
  \mbox{ then }
  P_{ij}^{\pm}(z,w) = 1
  \mbox{ and }
  d_{ij}=r.
\nonumber\\
&& \mbox{If }
  A_{i,\sigma(i)} = 0
  \mbox{ and }
  \sigma(j) = j,
  \mbox{ then }
  P_{ij}^{\pm}(z,w) =
  \frac
     {z^r q^{\pm 2r}-w^r}
     { z q^{\pm 2} -w}
  \mbox{ and }
  d_{ij} = r.
\nonumber\\
&& \mbox{If }
  A_{i,\sigma(i)} = 0
  \mbox{ and }
  \sigma(j) \neq j,
  \mbox{ then }
  P_{ij}^{\pm}(z,w) = 1
  \mbox{ and }
  d_{ij} = 1/2.
\nonumber\\
&& \mbox{If }
  A_{i,\sigma(i)}= -1,
  \mbox{ then }
  P_{ij}^{\pm}(z,w) =
  zq^{\pm r/2} + w
  \mbox{ and }
  d_{ij} = r/4. \nonumber
\end{eqnarray}
\end{theorem}

$\hfil\Box$

We denote by $\textstyle {\bf U}'$ the subalgebra of {\bf U} generated by
the elements $x_i^{\pm}(k)$, $a_i(l)$, $K_i^{\pm 1}$,
$\gamma^{\pm 1/2}$ $(i=1,2,\cdots,N, k\in {\bf Z},
l\in {\bf Z} \setminus \{0\})$.

We also need the following explicit isomorphism between the two
definitions.

\begin{prop} \cite{kn:D, kn:J3}\label{P:2.2}
 The ${\bf C}(q^{1/2})$-algebra isomorphism
$\Psi: U_q(X_N^{(r)}) \to \textstyle {\bf U}$ is given by
\begin{eqnarray}
&&e_i \mapsto x_i^{+}(0), \ \ f_i \mapsto x_i^{-}(0), \ \
t_i \mapsto K_i \ \ \mbox{if} \ \sigma(i)\neq i,
\nonumber\\
&&e_i \mapsto x_i^{+}(0), \ \ f_i \mapsto \frac 1r x_i^{-}(0), \ \
t_i \mapsto K_i \ \ \mbox{if} \ \sigma(i)=i,
\nonumber\\
 &&t_0 \mapsto \gamma K_{\theta}^{-1}, \ \ q^d \mapsto q^d,
\nonumber\\
&&e_0 \mapsto \frac 12[x_2^{-}(0), \cdots [x_{n-1}^{-}(0),
[x_{n}^{-}(0), [x_{n-1}^-(0), \cdots \nonumber\\
&&\mbox{}\hskip 1in \cdots[x_2^-(0), x_1^-(1)]_{q^{-1}} \cdots
]_{q^{-1}}]_{q^{-2}}]_{q^{-1}}\cdots ]_{q^{-1}} K_{\theta}^{-1},
\label{E:2.13}\\
&&f_0 \mapsto -q^{2n-2}K_{\theta}[x_2^{+}(0), \cdots [x_{n-1}^{+}(0),
[x_{n}^{+}(0), [x_{n-1}^+(0), \cdots \nonumber\\
&&\mbox{}\hskip 1in \cdots[x_2^+(0), x_1^+(-1)]_{q^{-1}} \cdots
]_{q^{-1}}]_{q^{-2}}]_{q^{-1}}\cdots ]_{q^{-1}}, \ \mbox{for }
A_{2n-1}^{(2)};\nonumber\\
&&e_0 \mapsto 2^{-n-1}[x_1^{-}(0), \cdots [x_{n-1}^{-}(0),
x_{n}^{-}(1)]_{q^{-2}} \cdots
]_{q^{-2}} K_{\theta}^{-1},
\nonumber\\
&&f_0 \mapsto (-1)^{n-1}q^{2n-2}K_{\theta}[x_1^{+}(0),  \cdots
[x_{n-1}^{+}(0),
x_{n}^{+}(-1)]_{q^{-2}} \cdots
]_{q^{-2}}, \ \mbox{for }
D_{n+1}^{(2)}; \nonumber\\
 &&e_0 \mapsto [2]^{-2n+2}[x_1^{-}(0), \cdots [x_n^{-}(0),
[x_{n}^{-}(0), [x_{n-1}^-(0), \cdots \nonumber\\
&&\mbox{}\hskip 1in \cdots[x_2^-(0),
x_1^-(1)]_{q^{-1}} \cdots
]_{q^{-1}}]_1]_{q^{-1}}\cdots ]_{q^{-1}}]_1 K_{\theta}^{-1},
\nonumber\\
&&f_0 \mapsto (-q)^{2n-3}[2]^{-2n+2}K_{\theta}
[x_1^{+}(0), \cdots [x_n^{+}(0),
[x_{n}^{+}(0), [x_{n-1}^-(0), \cdots \nonumber\\
&&\mbox{} \hskip 1in \cdots[x_{2}^+(0),
x_1^+(-1)]_{q^{-1}} \cdots
]_{q^{-1}}]_1]_{q^{-1}}\cdots ]_{q^{-1}}]_1,
\ \mbox{for }
A_{2n}^{(2)};\nonumber\\
&&e_0 \mapsto \frac 13[x_1^{-}(0),  [x_2^{-}(0),
x_1^-(1)]_{q^{-3}}]_{q^{-1}} K_{\theta}^{-1},
\nonumber\\
&&f_0 \mapsto q^4K_{\theta}[x_1^+(0), [x_{2}^{+}(0),
x_1^+(-1) ]_{q^{-3}}]_{q^{-1}}, \
\mbox{for } D_4^{(3)}.
\end{eqnarray}
where
$[x,y]_v=xy-vyx$ and ${\theta}$ is the maximal root
in the sense that $\delta=\theta+\a_0$ (in our convention):
$$
\theta=\left\{
     \begin{array}{ll}
      \a_1+2(\a_2+\cdots+\a_{n-1})+\a_n ,
                           &\mbox{for } A_{2n-1}^{(2)},\\
      \a_1+\cdots+\a_n,
                           &\mbox{for } D_{n+1}^{(2)},\\
      2(\a_1+\cdots+\a_n),
                           &\mbox{for } A_{2n}^{(2)},\\
      \a_1+2\a_2,          &\mbox{for } D_4^{(3)}.
     \end{array}
\right.
$$

The restriction of $\Psi$ to $U_q(X_N^{(r)})$ defines an isomorphism of
$U_q(X_N^{(r)})$ and ${\bf U}'$.  \hskip 1cm $\Box$
\end{prop}

Through the isomorphism of Proposition \ref{P:2.2} the coproduct will be
carried over
to the Drinfeld realization.

\begin{theorem}\label{T:2.3}
Let $k\in  {\bf Z}_{\ge 0}$, $l\in  {\bf N}$,
and let $N_{+}^s$ (resp. $N_{-}^s$) be the subalgebra of
{\bf U} generated by the elements $x_{i_1}^{+}(-m_1)\cdots x_{i_s}^{+}(-m_s)$
(resp. $x_{i_1}^{-}(m_1)\cdots x_{i_s}^{-}(m_s)$) with $m_i \in \textstyle
{\bf Z}_{\ge 0}$. Then the comultiplication $\Delta$ of the algebra {\bf U}
has the following form:
\begin{eqnarray*}
\Delta(x_{i}^{+}(k))&=&x_i^{+}(k) \otimes \gamma^k
+\gamma^{2k} K_i \otimes x_i^{+}(k)\\
&&\ \ \ \ \ \ \ \ + \sum_{j=0}^{k-1} \gamma^{\frac {k-j}{2}} \psi(k-j)
\otimes \gamma^{k-j} x_i^{+}(j) \ \ (\textstyle {mod} \ N_{-}\otimes N_{+}^2),
\\
 \Delta(x_{i}^{+}(-l))&=&x_i^{+}(-l) \otimes \gamma^{-l}
+ K_i^{-1} \otimes x_i^{+}(-l)\\
&&\ \ \ \ \ \ \ \ + \sum_{j=1}^{l-1} \gamma^{\frac {l-j}{2}} \varphi(-l+j)
\otimes \gamma^{-l+j} x_i^{+}(-j) \ \ (\textstyle {mod} \ N_{-}\otimes
N_{+}^2), \\
 \Delta(x_{i}^{-}(l))&=&x_i^{-}(l) \otimes K_i
+\gamma^{l} \otimes x_i^{-}(l)\\
&&\ \ \ \ \ \ \ \ + \sum_{j=1}^{l-1} \gamma^{l-j} x_i^{-}(j)
\otimes \gamma^{\frac {j-l}{2}} \psi_i(l-j)
 \ \ (\textstyle {mod} \ N_{-}^2 \otimes N_{+}), \\
 \Delta(x_{i}^{-}(-k))&=&x_i^{-}(k) \otimes \gamma^{-2k}K_i^{-1}
+\gamma^{-k} \otimes x_i^{-}(-k)\\
&&\ \ \ \ \ \ \ \ + \sum_{j=0}^{k-1} \gamma^{j-k} x_i^{-}(-j)
\otimes \gamma^{\frac {j-k}{2}} \varphi_i(j-k)
 \ \ (\textstyle {mod} \ N_{-}^2 \otimes N_{+}), \\
 \Delta(a_i(l))&=&a_i(l) \otimes \gamma^{\frac {l}{2}}
+ \gamma^{\frac {3l}{2}} \otimes a_i(l) \ \ (\textstyle {mod} \ N_{-}\otimes
N_{+}),\\
 \Delta(a_i(-l))&=&a_i(-l) \otimes \gamma^{-\frac {3l}{2}}
+ \gamma^{-\frac {l}{2}} \otimes a_i(-l) \ \ (\textstyle {mod} \ N_{-}\otimes
N_{+}).
\end{eqnarray*}
In the above formulas, we only take the nonzero generators:
$x_i^{\pm}(k)=a_i(k)=0$ if $\sigma(i)=i$ and
$k\neq 0\  mod\ r$.
Moreover the same formulas are true for the derived subalgebra
$\textstyle {\bf U}'$.
\end{theorem}

{\it Proof.} The proof is an induction on degrees similar to simply-laced
types \cite{kn:JKK}. The formulas of $\Delta(x^{\pm}_i(0))$ are clearly
valid as given
by the isomorphism in Proposition \ref{P:2.2} and the coproduct (\ref{E:2.6}).
To get the formulas for degree $\pm 1$ elements we use the inverse isomorphism
$\Psi^{-1}$: For a sequence $i=i_1, i_2, \cdots, i_{h^{(r)}-1}$, we have
\begin{equation}\label{E:2.15}
\begin{array}{rcl}
\Psi^{-1}(x_i^{+}(-1))&=& const \cdot t_{i_{2}}
[ f_{i_2}, t_{i_{3}}[f_{i_{3}},
\cdots t_{i_{h-1}}[f_{i_{h-1}}, t_{\theta}^{-1} f_0] \cdots ]], \\
\Psi^{-1}(x_i^{-}(1))&=& const\cdot
t_{i_{2}}^{-1} [ e_{i_{2}},
t_{i_{3}}^{-1} [e_{i_{3}},
\cdots t_{i_{h-1}}^{-1} [e_{i_{h-1}}, e_0 t_{\theta}] \cdots ]].
\end{array}
\end{equation}
Applying $\Delta$ it follows that
\begin{eqnarray*}
 \Delta(x_i^{+}(-1))&=&x_i^{+}(-1) \otimes \gamma^{-1}
+K_i^{-1} \otimes x_i^{+}(-1) \ \ (\textstyle {mod} \ N_{-} \otimes N_{+}^2),
\\
 \Delta(x_i^{-}(1))&=&x_i^{-}(1) \otimes K_i
+\gamma \otimes x_i^{-}(1) \ \ (\textstyle {mod} \ N_{-}^2 \otimes N_{+}).
\end{eqnarray*}
Using the relations
\begin{eqnarray*}
\ [x_{i}^{+}(0), x_{i}^{-}(1)]&=&p_i(q_i-q_i^{-1})^{-1} \gamma^{-1/2}
\psi_{i}(1)
=p_i\gamma^{-1/2}K_i a_i(1), \\
\ [x_{i}^{+}(-1), x_{i}^{-}(0)]&=&-p_i(q_i-q_i^{-1})^{-1} \gamma^{1/2}
\varphi_{i}(-1)
=p_i\gamma^{1/2}K_i^{-1} a_i(-1),
\end{eqnarray*}
where $p_i=1$ for $\sigma(i)\neq i$ and $p_i=r$ for $\sigma(i)=i$.
we obtain
\begin{eqnarray*}
 \Delta(a_i(1))&=&a_i(1) \otimes \gamma^{1/2} + \gamma^{3/2} \otimes a_i(1)
\ \ (\textstyle {mod} \ N_{-}\otimes N_{+}), \\
 \Delta(a_i(-1))&=&a_i(-1) \otimes \gamma^{-3/2} + \gamma^{-1/2}
\otimes a_i(-1) \ \ (\textstyle {mod} \ N_{-}\otimes N_{+}).
\end{eqnarray*}
Based on these two formulas and the Drinfeld relations it is easy to see by
induction that the formulas of $\Delta(x^{\pm}_i(k))$ are valid. Subsequently
the formulas for $\Delta(a_i(k))$ follow from the formulas for
$\Delta(\phi_i(k))$ and $\Delta(\psi_i(k))$, which proves the theorem.
\hskip 1cm $\Box$
\medskip

We now turn to some special finite dimensional
representations $V(X_N^{(r)})$ of
$U_q'(X_N^{(r)})$ arising from the level one
perfect crystal graphs. We list the representations case by case.
Note that in the case of $D_4^{(3)}$ the representation graph is different
from its crystal one.

For $A_{2n-1}^{(2)}$,
$$V=\left( \bigoplus_{i=1}^n {\bf C}(q^{1/2}) v_i \right) \oplus
\left( \bigoplus_{i=1}^n {\bf C}(q^{1/2}) v_{\overline i} \right).$$
The action
is given by:
\begin{equation}\label{E:A}
\begin{array}{rcl}
 e_i&=&E_{i, i+1}+E_{\overline {i+1}, \overline i}, \ \
f_{i}=E_{i+1, i}+E_{\overline i, \overline {i+1}}, \\
 t_i&=&q(E_{ii}+E_{\overline {i+1}, \overline {i+1}})
+q^{-1} (E_{\overline i, \overline i}+E_{i+1, i+1})
+\sum_{j\neq i, i+1, \overline i, \overline {i+1}} E_{jj},\\
 f_n&=&E_{\overline n, n}, \ \
e_n=E_{n, \overline n},\\
 t_n&=&q_nE_{nn}+q_n^{-1}E_{\overline {n}, \overline {n}}
+\sum_{j\neq n, \overline n} E_{jj} ,\\
f_0&=&E_{2\overline 1}+E_{1\overline 2}, \ \
e_0=E_{\overline 12}+E_{\overline 21}, \\
t_0&=&q(E_{\overline 1\overline 1}+E_{\overline 2\overline 2})
+q^{-1}(E_{1, 1}+E_{2,2})
+\sum_{j, \overline j\neq 1,2}E_{jj},
\end{array}
\end{equation}
where $i=1,2,\cdots, n-1$, $(d_0, \cdots, d_n)=
(1, \cdots, 1, 2)$, and $E_{ij}\in \textstyle{End}(V)$ so that
$
E_{ij}v_k=\delta_{jk}v_i
$.

For $D_{n+1}^{(2)}$,
$$V=\left( \bigoplus_{i=0}^n {\bf C}(q^{1/2}) v_i \right) \oplus
\left( \bigoplus_{i=0}^n {\bf C}(q^{1/2}) v_{\overline i} \right).$$
The action
is given by:
\begin{equation}\label{E:D}
\begin{array}{rcl}
 e_i&=&E_{i, i+1}+E_{\overline {i+1}, \overline i}, \ \
f_{i}=E_{i+1, i}+E_{\overline i, \overline {i+1}}, \\
 t_i&=&q_i(E_{ii}+E_{\overline {i+1}, \overline {i+1}})
+q_i^{-1} (E_{\overline i, \overline i}+E_{i+1, i+1})
+\sum_{j\neq i, i+1, \overline i, \overline {i+1}} E_{jj},\\
 e_n&=&[2]E_{n0}+E_{0\overline n}, \ \
f_n=[2]E_{\overline n0}+E_{0n},\\
 t_n&=&q^2E_{nn}+q^{-2}E_{\overline {n}, \overline {n}}
+\sum_{j\neq n, \overline n} E_{jj} ,\\
e_0&=&[2]E_{\overline 1\overline 0}+E_{\overline 01}, \ \
f_0=[2]E_{1\overline 0}+E_{\overline 0\overline 1}, \\
t_0&=&q^2E_{\overline 1\overline 1}+q^{-2}E_{11}
+\sum_{j\neq 1, \overline 1}E_{jj},
\end{array}
\end{equation}
where $(d_0, \cdots, d_n)=(1, 2, \cdots, 2, 1)$.

For $A_{2n}^{(2)}$,
$$V=\left( \bigoplus_{i=1}^n {\bf C}(q^{1/2}) v_i \right) \oplus
{\bf C}(q^{1/2})v_0\oplus
\left( \bigoplus_{i=1}^n {\bf C}(q^{1/2}) v_{\overline i} \right).$$
The action
is given by:
\begin{equation}\label{E:A_2}
\begin{array}{rcl}
 e_i&=&E_{i, i+1}+E_{\overline {i+1}, \overline i}, \ \
f_{i}=E_{i+1, i}+E_{\overline i, \overline {i+1}}, \\
 t_i&=&q_i(E_{ii}+E_{\overline {i+1}, \overline {i+1}})
+q_i^{-1} (E_{\overline i, \overline i}+E_{i+1, i+1})
+\sum_{j\neq i, i+1, \overline i, \overline {i+1}} E_{jj},\\
 e_n&=&[2]_nE_{n0}+E_{0\overline n}, \ \
f_n=[2]_nE_{\overline n0}+E_{0n},\\
 t_n&=&q_n^2E_{nn}+q_n^{-2}E_{\overline {n}, \overline {n}}
+\sum_{j\neq n, \overline n} E_{jj} ,\\
e_0&=&E_{\overline 11},  \ \
f_0=E_{1\overline 1}, \\
t_0&=&q_0(E_{\overline 1\overline 1})+q_0^{-1}(E_{11})
+\sum_{j\neq 1, \overline 1}E_{jj},
\end{array}
\end{equation}
where $(d_0, d_1, \cdots, d_n)=(2, 1, \cdots, 1, 1/2)$.
If $n=1$, then $(d_0, d_1)=(4, 1)$.

For $D_4^{(3)}$,
$$V=\left( \bigoplus_{i=0}^3 {\bf C}(q^{1/2}) v_i \right) \oplus
\left( \bigoplus_{i=0}^3 {\bf C}(q^{1/2}) v_{\overline i} \right).$$
The action
is given by:
\begin{equation}\label{E:D_4}
\begin{array}{rcl}
e_{1}&=&E_{12}+E_{\overline 2, \overline {1}}+E_{0\overline 3}+
[2]E_{30}+E_{3\overline 0}, \ \
f_{1}=E_{21}+E_{\overline 1, \overline {2}}+
E_{03}+[2]_1E_{\overline 30}+E_{\overline 3,\overline 0}, \\
 t_1&=&q(E_{11}+E_{\overline {2}\overline {2}})
+q^{-1} (E_{\overline 1\overline 1}+E_{22})+q^2E_{33}
+q^{-2}E_{\overline 3\overline 3}+E_{00}+E_{\overline 0\overline 0}\\
e_2&=&E_{23}+E_{\overline 3\overline 2}, \ \
f_2=E_{32}+E_{\overline 2\overline 3},\\
 t_2&=&q_2(E_{22}+E_{\overline 3\overline 3})+
 q_2^{-1}(E_{33}+E_{\overline 2\overline 2})
+\sum_{j=0,1, \overline 0, \overline 1} E_{jj} ,\\
e_0&=&[2]E_{\overline 1\overline 0}+E_{\overline 01}+
   E_{\overline 32}+E_{\overline 23}+E_{\overline 1, 0}, \ \
f_0=[2]E_{1\overline 0}+E_{\overline 0\overline 1}+
E_{2\overline 3}+E_{3\overline 2}+E_{10}, \\
t_0&=&q^{-1}(E_{22}+E_{33})+q(E_{\overline 2\overline 2}+
E_{\overline 3\overline 3})+q^{-2}E_{11}+q^2E_{\overline 1\overline 1}
+E_{00}+E_{\overline 0\overline 0},
\end{array}
\end{equation}
where $(d_0, d_1, d_2)=(1, 1, 3)$. Note that the associated
crystal with the last matrices
of $e_1, f_1, e_0, f_0$ removed, is a perfect crystal for $D_4^{(3)}$.
The perfect crystal structure (see \cite{kn:KMN}) is more than we will
need in the sequel.

Since $V$ is a finite dimensional vector space over ${\bf C}
(q^{1/2})$, it does not admit a $\U$-module structure. However
we can
define a $\U$-module structure on the affinization or evaluation
module of $V$.
The affinization of $V$ is the $\U$-module
$V_{z}=V \otimes {\bf C}(q^{1/2})[z, z^{-1}]$ with the following
action:
\begin{eqnarray}
 e_i(v\otimes z^m)&=&e_i v \otimes z^{m+\delta_{i,0}}, \ \
 f_i(v\otimes z^m)=f_i v \otimes z^{m-\delta_{i,0}},
 \nonumber\\
 t_i(v\otimes z^m)&=&t_i v \otimes z^{m},  \label{E:2.18}\\
 q^d (v\otimes z^m)&=& q^m v\otimes z^m,   \nonumber
\end{eqnarray}
for $i=0,1,\cdots,n$, $m\in {\bf Z}$, $v\in V$.

The evaluation module $V_z$ is a level zero $\U$-module, i.e., $\gamma$ acts
as identity ($=q^0$). Through the isomorphism $\Psi$ the evaluation module
is also a $\bf U$-module. The action of the Drinfeld generators are given
by the following result.

\begin{theorem}\label{T:2.4} The Drinfeld generators act on the
evaluation module $V_z$ as follows.

For $V(A_{2n-1}^{(2)})$:
\begin{equation}
\begin{array}{rcl}
x_i^{+}(k)&=&(q^iz)^k E_{i,i+1}+(-q^{2n-i-1}z)^k E_{\overline {i+1},
\overline i},\\
x_i^{-}(k)&=&(q^iz)^k E_{i+1,i}+(-q^{2n-i-1}z)^k E_{\overline i,
\overline {i+1}},\\
x_n^{+}(2k)&=&(q^{n}z)^{2k}E_{n,\overline n}, \\
x_n^{-}(2k)&=&2(q^{n}z)^{2k}E_{\overline n, n}, \\
a_i(l)&=&\displaystyle\frac {[l]}{l} ( (q^iz)^l (q^{-l}E_{i, i}
-q^l E_{i+1,i+1})\\
&& \ \ \ \ \ \ \ \  +(-q^{2n-i}z)^l (q^{-l}E_{\overline {i+1},
\overline {i+1}}
-q^l E_{\overline i, \overline i}) ),\\
a_n(2l)&=&\displaystyle\frac {[l]_n}{l}  (q^{n}z)^{2l} ( q^{-2l}E_{n,n}
-q^{2l} E_{\overline {n},\overline {n}}) ,
\end{array}
\end{equation}
where $i=1, 2, \cdots, n-1$, $k\in {\bf Z}$, and $l\in {\bf Z}
\setminus \{0\}$, $x^{\pm}_n(m)=a_n(m)=0$ for $m$ odd, and $q_n=q^2$.

For $V(D_{n+1}^{(2)})$:
\begin{equation}
\begin{array}{rcl}
x_i^{+}(2k)&=&(-q^2)^{kn}z^{2k}
(q^{-2(n-i)k}E_{i,i+1}+q^{2(n-i)k} E_{\overline {i+1},
\overline i}), \\
x_i^{-}(2k)&=&2(-q^2)^{kn}z^{2k}
(q^{-2(n-i)k}E_{i+1,i}+q^{2(n-i)k} E_{\overline {i},
\overline {i+1}}), \\
a_i(2l)&=&\displaystyle\frac {[l]_i}{l} (-q^2)^{ln}z^{2l}
( q^{-2(n-i+1)l}E_{i,i}-q^{-2(n-i-1)l} E_{i+1,i+1}\\
&& \ \ \ \ \ \ \ \  +q^{2(n-i-1)l}E_{\overline {i+1},
\overline {i+1}}
-q^{2(n-i+1)l} E_{\overline i, \overline i} ),\\
x_n^{+}(2k)&=&(-q^2)^{kn}z^{2k}([2]E_{n,0}+E_{0,\overline n}), \\
x_n^{-}(2k)&=&(-q^2)^{kn}z^{2k}(E_{0,n}+[2]E_{\overline n,0}), \\
a_n(2l)&=&(-q^2)^{nl}z^{2l}
\displaystyle\frac {[2l]}{2l} ( 2q^{-2l}E_{n,n}
-2q^{2l} E_{\overline {n},\overline {n}}+(q^{-2l}-q^{2l})(E_{0,0}-
E_{\overline 0, \overline 0}) ) ,\\
x_n^{+}(2k+1)&=&(-q^2)^{kn}z^{2k+1}((-q^2)^nE_{\overline 0,\overline n}
-[2]E_{n, \overline 0}), \\
x_n^{-}(2k+1)&=&(-q^2)^{kn}z^{2k+1}([2]E_{\overline n,\overline 0}
-(-q^2)^nE_{\overline 0, n}),  \\
a_n(2l+1)&=&(-q^2)^{nl}z^{2l+1}
\displaystyle\frac {[4l+2]}{2l+1} ( E_{0,\overline {0}}+
(-q^2)^nE_{\overline 0, 0} )
\end{array}
\end{equation}
where $i=1, 2, \cdots, n-1$, $k\in {\bf Z}$, and $l\in {\bf Z}
\setminus \{0\}$, $x^{\pm}_i(m)=a_i(m)=0$ for $m$ odd, and $q_i=q^2$.

For $V(A_{2n}^{(2)})$:
\begin{equation}
\begin{array}{rcl}
x_i^{+}(k)&=&(q^iz)^k (E_{i,i+1}+(-q^{2n-2i-1}z)^k E_{\overline {i+1},
\overline i}),\\
x_i^{-}(k)&=&(q^iz)^k (E_{i+1,i}+(-q^{2n-2i-1}z)^k E_{\overline i,
\overline {i+1}}),\\
a_i(l)&=&\displaystyle\frac {[l]}{l}  (q^iz)^l (q^{-l}E_{i, i}
-q^l E_{i+1,i+1}\\
&& \ \ \ \ \ \ \ \  +(q^{2n-2i+1}z)^l (q^{-l}E_{\overline {i+1},
\overline {i+1}}
-q^l E_{\overline i, \overline i}) ),\\
x_n^{+}(k)&=&(q^{n}z)^{k}((-q)^kE_{0,\overline n}+[2]_nE_{n,0}), \\
x_n^{-}(k)&=&(q^{n}z)^{k}(E_{0,n}+(-q)^k[2]_nE_{\overline n, 0}, \\
a_n(l)&=&\displaystyle\frac {[2l]_n}{l}  (q^{n}z)^{l} ( q^{-l}E_{n,n}
-q^{l}E_{0,0}+(-q)^l(q^{-l}E_{0,0}-q^lE_{\overline {n},\overline {n}}) ) ,
\end{array}
\end{equation}
where $i=1, 2, \cdots, n-1$, $k\in {\bf Z}$, and $l\in {\bf Z}
\setminus \{0\}$, $q_0=q^2$, $q_i=q$ and $q_n=q^{1/2}$.

For $V(D_4^{(3)})$:
\begin{equation}
\begin{array}{rcl}
x_1^{+}(3k)&=&(q^iz)^{3k}(E_{12}+q^{12k} E_{\overline {2}
\overline 1}+q^{6k}E_{0\overline 3}+q^{6k}E_{3\overline 0}+q^{6k}[2]
E_{30}),\\
x_1^{+}(3k+1)&=&(q^iz)^{3k+1}(E_{12}+q^{12k+4} E_{\overline {2}
\overline 1}+q^{6k+2}(q^2-1)E_{0\overline 3}-q^{6k+5}E_{\overline{03}}\\
\ && \ \ \ \ -q^{6k}(q^2+q^{-2})E_{3\overline 0}-q^{6k+3}
E_{30}),\\
x_1^{+}(3k+2)&=&(q^iz)^{3k}(E_{12}+q^{12k+8} E_{\overline {2}
\overline 1}-q^{6k+6}E_{0\overline 3}+q^{6k+7}E_{\overline{03}}+
q^{6k}E_{3\overline 0}+q^{6k+3}
E_{30}),\\
x_1^{-}(3k)&=&(q^iz)^{3k}(E_{21}+q^{12k} E_{\overline {1}
\overline 2}+q^{6k}E_{\overline{30}}+q^{6k}[2]E_{\overline 30}+q^{6k}
E_{03}),\\
x_1^{-}(3k+1)&=&(q^iz)^{3k+1}(E_{21}+q^{12k+4} E_{\overline {1}
\overline 2}+q^{6k-2}E_{\overline{30}}\\
\ &&\ \ \ \ -q^{6k+1}E_{\overline 30}+
q^{6k+5}E_{\overline 03}-q^{6k+4}
E_{03}),\\
x_1^{-}(3k+2)&=&(q^iz)^{3k+2}(E_{21}+q^{12k+8} E_{\overline {1}
\overline 2}-q^{6k+2}(q^2+q^{-2})E_{\overline{30}}-q^{6k+5}E_{\overline 30}\\
\ && \ \ \ \ -q^{6k+7}E_{\overline 03}+q^{6k+4}(q^2-1)
E_{03}),\\
x_2^{+}(3k)&=&(q^2z)^{3k}(E_{23}+q^{6k}E_{\overline{32}}), \\
x_2^{-}(3k)&=&3(-1)^k(q^{2}z)^{3k}(E_{32}+q^{6k}E_{\overline{23}}), \\
a_1(3l)&=&\displaystyle\frac {[3l]}{3l} (qz)^{3l}( q^{-3l}E_{11}-q^{15l}
E_{\overline{11}}-q^{3l}E_{22}+q^{9l}E_{\overline{22}}+2q^{3l}E_{33}\\
\ &&\ \ \ \ +(q^{3l}-q^{9l})E_{00}+(q^{3l}-q^{9l})E_{\overline{00}})\\
a_1(3k+1)&=&\displaystyle\frac{[3k+1]}{3k+1}(qz)^{3k+1}(q^{-3k-1}E_{11}
-q^{15k+5}E_{\overline{11}}-q^{3k+1}E_{22}+q^{9k+3}E_{\overline{22}}
\\
\ &&\ \ \ \ +q^{9k+3}E_{\overline{33}}
-q^{3k+1}E_{33}+\frac{q^3+q^{-3}}{q+q^{-1}}(q^{3k+1}+q^{-3k-1})
q^{6k+1}E_{0\overline 0}\\
\ &&\ \ \ \ -(q^{3k+1}+q^{-3k-1})q^{6k+5}E_{\overline 00}
 +(q^{9k+5}+q^{3k+3}-q^{3k+1})E_{00}\\
\ &&\ \ \ \ +
(-q^{9k+5}-q^{3k+3}+q^{9k+3})E_{\overline{00}})\\
a_1(3k+2)&=&\displaystyle\frac{[3k+2]}{3k+2}(qz)^{3k+2}(q^{-3k-2}E_{11}
-q^{15k+10}E_{\overline{11}}-q^{3k+2}E_{22}+q^{9k+6}E_{\overline{22}}\\
\ &&\ \ \ \ +q^{9k+6}E_{\overline{33}}-q^{3k+2}E_{33}-
\frac{q^3+q^{-3}}{q+q^{-1}}(q^{3k+2}+q^{-3k-2})
q^{6k+3}E_{0\overline 0}\\
\ &&\ \ \ \ +(q^{3k+2}+q^{-3k-2})q^{6k+7}E_{\overline 00}
+(q^{9k+6}+q^{9k+8}-q^{3k+4})E_{00}\\
\ && \ \ \ \ +
(q^{9k+8}-q^{3k+4}+q^{3k+2})E_{\overline{00}})\\
a_2(3l)&=&\displaystyle\frac {[l]_2}{l}  (q^2z)^{3l} ( q^{-3l}E_{22}
-q^{3l}E_{33}+q^{6l} (q^{-3l}E_{\overline{33}}-q^{3l}E_{\overline{22}}) ,
\end{array}
\end{equation}
where $k\in {\bf Z}$, and $l\in {\bf Z}
\setminus \{0\}$, $q_0=q_1=q$ and $q_2=q^3$.
\end{theorem}

{\it Proof.} The theorem is proved by induction using the isomorphism of
two definitions. Note that the formulas of degree zero hold by
the definitions (\ref{E:A}), (\ref{E:D}), (\ref{E:A_2}), and
(\ref{E:D_4}).
Now using (\ref{E:2.15})
we can get degree one relations right, then higher degree relations follow
using Drinfeld relations.
\hskip 1cm $\Box$

\begin{theorem}\label{T:2.5} The Drinfeld generators act on the
dual evaluation module $V_z^*=V^*\otimes {\bf C}(q^{1/2})[z, z^{-1}]$
as follows.

For $V(A_{2n-1}^{(2)})$:
\begin{equation}\label{E:2.21}
\begin{array}{rcl}
x_i^{+}(k)&=&(-q^{-1})((q^{-i}z)^k E_{i+1,i}^*+(-q^{-2n+i}z)^k
E_{\overline {i},
\overline {i+1}}^*),\\
x_i^{-}(k)&=&(-q)((q^{-i}z)^k E_{i,i+1}^*+(-q^{-2n+i}z)^k E_{\overline {i+1},
\overline {i}}^*),\\
a_i(l)&=&\displaystyle\frac {[l]}{l} ((q^{-i}z)^l(q^{-l} E^*_{i+1,i+1}
- q^{l}E_{ii}^*)\\
&& \ \ \ \ \ \ \ \  +(-q^{-2n+i}z)^l(q^{-l} E_{\overline i, \overline i}^*
- q^{l}E_{\overline {i+1},
\overline {i+1}}^* )),\\
x_n^{+}(2k)&=&(-q^{-2})((q^{-n}z)^{2k}E_{\overline n, n}^*,\\
x_n^{-}(2k)&=&2(-q^2)((q^{-n}z)^{2k}E_{n, \overline n}^*,\\
a_n(2l)&=&\displaystyle\frac {[l]_n}{l} (q^{-n}z)^{2l}
((q^{-2l} E_{\overline {n}, \overline {n}}^*-q^{2l}E_{n, n}^*),
\end{array}
\end{equation}
where $i=1, 2, \cdots, n-1$, $k\in {\bf Z}$, and $l\in {\bf Z}
\setminus \{0\}$, and $x_n^{\pm}(2k+1)=a_n(2k+1)=0$.

For $V(D_{n+1}^{(2)})$:
\begin{equation}%
\begin{array}{rcl}
x_i^{+}(2k)&=&-q^{-2}(-q^2)^{-kn}z^{2k}
(q^{2(n-i)k}E_{i+1, i}^*+q^{-2(n-i)k} E_{\overline {i},
\overline {i+1}}^*), \\
x_i^{-}(2k)&=&-2q^2(-q^2)^{-kn}z^{2k}
(q^{2(n-i)k}E_{i,i+1}^*+q^{-2(n-i)k} E_{\overline {i+1},
\overline {i}}^*), \\
a_i(2l)&=&\displaystyle\frac {[l]_i}{l} (-q^2)^{-ln}z^{2l}
( q^{2(n-i-1)l}E_{i+1,i+1}^*-q^{2(n-i+1)l} E_{i,i}^*\\
&& \ \ \ \ \ \ \ \  +q^{-2(n-i+1)l}E_{\overline {i},
\overline {i}}
-q^{-2(n-i-1)l} E_{\overline {i+1}, \overline {i+1}} ),\\
x_n^{+}(2k)&=&(-q^2)^{-kn-1}z^{2k}([2]q^{-2}E_{0, n}^*+E_{\overline n, 0}^*),
\\
x_n^{-}(2k)&=&(-q^2)^{-kn+1}z^{2k}(q^2E_{n,0}^*+[2]E_{0, \overline n}^*), \\
a_n(2l)&=&(-q^2)^{-nl}z^{2l}
\displaystyle\frac {[2l]}{2l} ( 2q^{-2l}E_{\overline n,\overline n}^*
-2q^{2l} E_{n, n}^*+(q^{-2l}-q^{2l})(E_{0,0}^*-
E_{\overline 0, \overline 0}^*) ) ,\\
x_n^{+}(2k+1)&=&(-q^2)^{-kn-1}z^{2k+1}((-q^2)^{-n}
E_{\overline n,\overline 0}^*
-[2]q^{-2}E_{\overline 0, n}^*), \\
x_n^{-}(2k+1)&=&(-q^2)^{-kn+1}z^{2k+1}([2]E_{\overline 0,\overline n}^*
-q^2(-q^2)^{-n}E_{n, \overline 0}^*),  \\
a_n(2l+1)&=&-(-q^2)^{-nl}z^{2l+1}
\displaystyle\frac {[4l+2]}{2l+1} ( E_{\overline {0}, 0}^*+
(-q^2)^{-n}E_{0, \overline 0}^* )
\end{array}
\end{equation}
where $i=1, 2, \cdots, n-1$, $k\in {\bf Z}$, and $l\in {\bf Z}
\setminus \{0\}$, $x^{\pm}_i(m)=a_i(m)=0$ for $m$ odd, and $q_i=q^2$.

For $V(A_{2n}^{(2)})$:
\begin{equation}
\begin{array}{rcl}
x_i^{+}(k)&=&-q^{-1}(q^{-i}z)^k (E_{i+1,i}^*+
(-q^{-2n+2i+1}z)^k E_{\overline {i},
\overline {i+1}}^*),\\
x_i^{-}(k)&=&-q(q^{-i}z)^k (E_{i,i+1}^*+(-q^{-2n+2i+1}z)^k E_{\overline {i+1},
\overline {i}}^*),\\
a_i(l)&=&\displaystyle\frac {[l]}{l} ( (q^{-i}z)^l (q^{-l}E_{i+1, i+1}^*
-q^l E_{i,i}^*)\\
&& \ \ \ \ \ \ \ \  +(q^{-2n+2i-1}z)^l (q^{-l}E_{\overline {i},
\overline {i}}^*
-q^l E_{\overline {i+1}, \overline {i+1}}^*) ),\\
x_n^{+}(k)&=&-q^{-1}(q^{-n}z)^{k}((-q)^{-k}E_{\overline n,0}^*+
q^{-1}[2]_nE_{0,n}^*), \\
x_n^{-}(k)&=&-q(q^{-n}z)^{k}(qE_{n,0}^*+(-q)^{-k}[2]_nE_{0,\overline n}^*, \\
a_n(l)&=&\displaystyle\frac {[2l]_n}{l}  (q^{-n}z)^{l} ( q^{-l}(E_{0,0}^*
-q^{l}E_{n,n}^*+(-q)^{-l}(q^{-l}E_{\overline n, \overline n}^*-
q^lE_{0,0}) ) ,
\end{array}
\end{equation}
where $i=1, 2, \cdots, n-1$, $k\in {\bf Z}$, and $l\in {\bf Z}
\setminus \{0\}$, $q_0=q^2$, $q_i=q$ and $q_n=q^{1/2}$.

For $V(D_4^{(3)})$:
\begin{equation}%
\begin{array}{rcl}
x_1^{+}(3k)&=&(-q^{-1})(q^iz)^{3k}(E_{21}^*+q^{12k} E_{\overline 1
\overline 2}^*+q^{6k+1}E_{\overline 30}^*+
q^{6k-1}E_{\overline 03}^*+q^{6k-1}[2]
E_{03}^*),\\
x_1^{+}(3k+1)&=&(-q^{-1})(q^iz)^{3k+1}(E_{21}^*+q^{12k+4} E_{\overline 1
\overline 2}^*+q^{6k+2}(q^2-1)E_{\overline 30}^*-q^{6k+6}E_{3\overline 0}^*\\
\ && \ \ \ \ -q^{6k-1}(q^2+q^{-2})E_{\overline 03}^*-q^{6k+2}
E_{03}^*),\\
x_1^{+}(3k+2)&=&(-q^{-1})(q^iz)^{3k}(E_{21}^*+q^{12k+8} E_{\overline{12}}^*
-q^{6k+7}E_{\overline 30}^*+q^{6k+8}E_{\overline{30}}^*\\
\ &&\ \ \ \ +
q^{6k-1}E_{\overline 03}^*+q^{6k+2}
E_{03}^*),\\
x_1^{-}(3k)&=&(-q)(q^iz)^{3k}(E_{12}^*+q^{12k} E_{\overline {2}
\overline 1}^*+q^{6k-1}E_{\overline{03}}^*+q^{6k-1}[2]E_{0\overline 3}^*
+q^{6k+1}
E_{30}^*),\\
x_1^{-}(3k+1)&=&(-q)(q^iz)^{3k+1}(E_{12}^*+q^{12k+4} E_{\overline {2}
\overline 1}^*+q^{6k-3}E_{\overline{03}}^*\\
\ &&\ \ \ \ -q^{6k}E_{0\overline 3}^*+
q^{6k+6}E_{3\overline 0}^*-q^{6k+5}
E_{30}^*),\\
x_1^{-}(3k+2)&=&(-q)(q^iz)^{3k+2}(E_{12}^*+q^{12k+8} E_{\overline {2}
\overline 1}^*-q^{6k+1}(q^2+q^{-2})E_{\overline{03}}^*
-q^{6k+4}E_{0\overline 3}^*\\
\ && \ \ \ \ -q^{6k+8}E_{3\overline 0}^*+q^{6k+5}(q^2-1)
E_{30}^*),\\
x_2^{+}(3k)&=&(-q^{-3})(q^2z)^{3k}(E_{32}^*+q^{6k}E_{\overline{23}}^*), \\
x_2^{-}(3k)&=&3(-q^3)(-1)^k(q^{2}z)^{3k}(E_{23}^*+q^{6k}E_{\overline{32}}^*),
\\
a_1(3l)&=&\displaystyle\frac {[3l]}{3l} (q^{-1}z)^{3l}( -q^{3l}E_{11}^*
+q^{-15l}
E_{\overline{11}}^*+q^{-3l}E_{22}^*-q^{-9l}E_{\overline{22}}^*
-2q^{-3l}E_{33}^*\\
\ &&\ \ \ \ -(q^{-3l}-q^{-9l})E_{00}^*-(q^{-3l}-q^{-9l})E_{\overline{00}}^*)\\
a_1(3k+1)&=&\displaystyle\frac{[3k+1]}{3k+1}(q^{-1}z)^{3k+1}
(-q^{3k+1}E_{11}^*
+q^{-15k-5}E_{\overline{11}}^*+q^{-3k-1}E_{22}^*\\
\ &&\ \ \ \ -q^{-9k-3}E_{\overline{22}}^*-q^{-9k-3}E_{\overline{33}}^*
-q^{-3k-1}E_{33}^*\\
\ &&\ \ \ \ -\frac{q^3+q^{-3}}{q+q^{-1}}(q^{-3k-1}+q^{3k+1})
q^{-6k-1}E_{\overline 00}^*\\
\ &&\ \ \ \ +(q^{-3k-1}+q^{3k+1})q^{-6k-5}E_{0\overline 0}^*
 -(q^{-9k-5}+q^{-3k-3}-q^{-3k-1})E_{00}^*\\
\ &&\ \ \ \  -(-q^{-9k-5}-q^{-3k-3}+q^{-9k-3})E_{\overline{00}}^*)\\
a_1(3k+2)&=&\displaystyle\frac{[3k+2]}{3k+2}(q^{-1}z)^{3k+2}(-q^{3k+2}E_{11}^*
+q^{-15k-10}E_{\overline{11}}^*+q^{-3k-2}E_{22}^*
-q^{-9k-6}E_{\overline{22}}^*\\
\ &&\ \ \ \ -q^{-9k-6}E_{\overline{33}}^*+q^{-3k-2}E_{33}^*+
\frac{q^3+q^{-3}}{q+q^{-1}}(q^{-3k-2}+q^{3k+2})
q^{-6k-3}E_{\overline 00}^*\\
\ &&\ \ \ \ -(q^{-3k-2}+q^{3k+2})q^{-6k-7}E_{0\overline 0}^*
-(q^{-9k-6}+q^{-9k-8}-q^{-3k-4})E_{00}^*\\
\ && \ \ \ \ -
(q^{-9k-8}-q^{-3k-4}+q^{-3k-2})E_{\overline{00}}^*)\\
a_2(3l)&=&\displaystyle\frac {[l]_2}{l}  (q^{-2}z)^{3l} ( -q^{3l}E_{22}^*
+q^{-3l}E_{33}^*-q^{-6l} (q^{3l}E_{\overline{33}}^*
+q^{-3l}E_{\overline{22}}^*)) ,
\end{array}
\end{equation}

\end{theorem}
{\it Proof.} Recall that the dual module $W^*$ of a module $W$ is defined by
$$
a.u^*(v)=u^*(S(a).v)
$$
for $a\in \U$, $u^*\in W^*$ and $v\in W$.  Explicitly we have the following
description of the dual modules. We only write down the matrices for the
generators. The underlying space is the standard dual vector space $V^*$
spanned by $v_i^*$, $i$ runs through the index set of various space $V$.
We also use $E_{ij}^*\in End(V^*)$ so that $E_{ij}^*v_k^*=\delta_{jk}v_i^*$.

For $A_{2n-1}^{(2)}$,
$$V=\left( \bigoplus_{i=1}^n {\bf C}(q^{1/2}) v_i \right) \oplus
\left( \bigoplus_{i=1}^n {\bf C}(q^{1/2}) v_{\overline i} \right).$$
The action
is given by:
\begin{equation}
\begin{array}{rcl}
 e_i&=&-q^{-1}(E_{i+1, i}^*+E_{\overline {i}, \overline {i+1}}^*), \ \
f_{i}=-q(E_{i, i+1}^*+E_{\overline {i+1}, \overline {i}}^*), \\
 t_i&=&q^{-1}(E_{ii}^*+E_{\overline {i+1}, \overline {i+1}}^*)
+q (E_{\overline i, \overline i}^*+E_{i+1, i+1}^*)
+\sum_{j\neq i, i+1, \overline i, \overline {i+1}} E_{jj}^*,\\
 f_n&=&-q_n^{-1}E_{n, \overline n}^*, \ \
e_n=-q_nE_{\overline n, n}^*,\\
 t_n&=&q_n^{-1}E_{nn}^*+q_nE_{\overline {n}, \overline {n}}^*
+\sum_{j\neq n, \overline n} E_{jj}^* ,\\
e_0&=&-q^{-1}(E_{2\overline 1}^*+E_{1\overline 2}^*), \ \
f_0=-q(E_{\overline 12}^*+E_{\overline 21}^*), \\
t_0&=&q^{-1}(E_{\overline 1\overline 1}^*+E_{\overline 2, \overline 2}^*)
+q(E_{1,1}^*+E_{2,2}^*)
+\sum_{j, \overline j\neq 1,2}E_{jj}^*.
\end{array}
\end{equation}
where $i=1,2,\cdots, n-1$, $(d_0, \cdots, d_n)=
(1, \cdots, 1, 2)$.

For $D_{n+1}^{(2)}$,
\begin{equation}%
\begin{array}{rcl}
 e_i&=&-q^{-2}(E_{i+1, i}^*+E_{\overline {i}, \overline {i+1}}^*), \ \
f_{i}=-q^2(E_{i, i+1}^*+E_{\overline {i+1}, \overline {i}}^*, \\
 t_i&=&q^{-2}(E_{ii}^*+E_{\overline {i+1}, \overline {i+1}}^*)
+q^{2} (E_{\overline i, \overline i}^*+E_{i+1, i+1}^*)
+\sum_{j\neq i, i+1, \overline i, \overline {i+1}} E_{jj}^*,\\
 e_n&=&-q^{-2}([2]E_{0,n}^*+q^{2}E_{\overline n,0}^*), \ \
f_n=-q^2(q^{-2}[2]E_{0, \overline n}^*+E_{n,0}^*),\\
 t_n&=&q^{-2}E_{nn}^*+q^{2}E_{\overline {n}, \overline {n}}^*
+\sum_{j\neq n, \overline n} E_{jj}^* ,\\
e_0&=&-([2]q^{-2}E_{\overline 0\overline 1}^*+E_{1, \overline 0}^*), \ \
f_0=-([2]E_{\overline 0, 1}^*+q^2E_{\overline 1\overline 0}^*), \\
t_0&=&q^{-2}E_{\overline 1\overline 1}^*+q^{2}(E_{11}^*
+\sum_{j\neq 1, \overline 1}E_{jj}^*.
\end{array}
\end{equation}
where $(d_0, \cdots, d_n)=(1, 2, \cdots, 2, 1)$.

For $A_{2n}^{(2)}$,
\begin{equation}%
\begin{array}{rcl}
 e_i&=&-q^{-1}(E_{i+1, i}^*+E_{\overline {i}, \overline {i+1}}^*), \ \
f_{i}=-q(E_{i, i+1}^*+E_{\overline {i+1}, \overline {i}}^*), \\
 t_i&=&q^{-1}(E_{ii}^*+E_{\overline {i+1}, \overline {i+1}}^*)
+q (E_{\overline i, \overline i}^*+E_{i+1, i+1}^*)
+\sum_{j\neq i, i+1, \overline i, \overline {i+1}} E_{jj}^*,\\
 e_n&=&-q^{-1}([2]_nE_{0, n}^*+qE_{\overline n,0}^*), \ \
f_n=-q(q^{-1}[2]_nE_{0, \overline n}^*+E_{n,0}^*),\\
 t_n&=&q^{-1}E_{nn}^*+qE_{\overline {n}, \overline {n}}^*
+\sum_{j\neq n, \overline n} E_{jj}^* ,\\
e_0&=&-q^{-2}E_{1, \overline 1}^*,  \ \
f_0=-q^2E_{\overline 1, 1}^*, \\
t_0&=&q_0^{-1}E_{\overline 1\overline 1}^*+q_0E_{11}^*
+\sum_{j\neq 1, \overline 1}E_{jj}^*.
\end{array}
\end{equation}
where $(d_0, d_1, \cdots, d_n)=(2, 1, \cdots, 1, 1/2)$.

For $D_4^{(3)}$,
\begin{equation}%
\begin{array}{rcl}
e_{1}&=&(-q^{-1})(E_{21}^*+E_{\overline 1, \overline {2}}^*+
qE_{\overline 30}^*+q^{-1}[2]E_{03}^*+q^{-1}E_{\overline 03}^*), \\
f_{1}&=&-q(E_{12}^*+E_{\overline 2, \overline {1}}^*+
qE_{30}^*+q^{-1}[2]_1E_{0\overline 3}^*+q^{-1}E_{\overline{03}}^*), \\
 t_1&=&q^{-1}(E_{11}^*+E_{\overline {2}\overline {2}}^*)
+q (E_{\overline 1\overline 1}^*+E_{22}^*)+q^{-2}E_{33}^*
+q^{2}E_{\overline 3\overline 3}^*+E_{00}^*+E_{\overline 0\overline 0}^*\\
e_2&=&(-q^{-3})(E_{32}^*+E_{\overline 2\overline 3}^*), \ \
f_2=(-q^3)(E_{23}^*+E_{\overline 3\overline 2}^*),\\
 t_2&=&q^{-3}(E_{22}^*+E_{\overline 3\overline 3}^*)+
 q^3(E_{33}^*+E_{\overline 2\overline 2}^*)
+\sum_{j=0,1, \overline 0, \overline 1} E_{jj}^* ,\\
e_0&=&(-q^{-1}z)(q^{-1}[2]E_{\overline 0\overline 1}^*+qE_{1\overline 0}^*+
   E_{2\overline 3}^*+E_{3\overline 2}^*+q^{-1}E_{0\overline 1}^*), \\
f_0&=&(-qz^{-1})(q^{-1}[2]E_{\overline 01}^*+qE_{\overline 1\overline 0}+
E_{\overline 32}^*+E_{\overline 23}^*+q^{-1}E_{01}^*), \\
t_0&=&q(E_{22}^*+E_{33}^*)+q^{-1}(E_{\overline 2\overline 2}^*+
E_{\overline 3\overline 3}^*)+q^{2}E_{11}^*+q^{-2}E_{\overline 1\overline 1}^*
+E_{00}^*+E_{\overline 0\overline 0}^*.
\end{array}
\end{equation}
where $(d_0, d_1, d_2)=(1, 1, 3)$.

In the above formulas, $i=1,2,\cdots, n-1$. They show that the Drinfeld
generators of degree zero elements do act on $V^*$ as given in the theorem.
The higher degree generators are obtained by inductions using Drinfeld
relations by the same argument as in Theorem \ref{T:2.4}.
\hskip 1cm $\Box$

\section{Level one realizations of $U_q(X^{(r)}_N)$}

In this section we recall the construction of level one twisted quantum
affine algebras given in \cite{kn:J1}. The level one integrable
modules are: $V(\Lambda_0), V(\Lambda_1)$ for $A_{2n-1}^{(2)}$;
$V(\Lambda_0), V(\Lambda_n)$ for $D_{n+1}^{(2)}$; $V(\Lambda_n)$ for
$A_{2n}^{(2)}$; and $V(\Lambda_0)$ for $D_4^{(3)}$.

Let $Q'$ be the root lattice of the simply laced finite dimensional Lie
algebra of types $A_{2n-1}, D_{n+1}, A_{2n}, D_4$. We shall use
the same convention or notation as in section one.
The diagram automorphism $\sigma$ decomposes
$Q'={\bf Z}\alpha_1'\oplus\cdots \oplus{\bf Z}\alpha_N'$
into orbits: $Q'=Q_0'\oplus Q_1'\cdots \oplus Q_{r-1}'$. We shall also
use $Q$ to denote the fixed-point sublattice $Q'_0$. We extend
the action of $\sigma$ to the weight lattice $P'={\bf Z}\la_1'\oplus\cdots
\oplus{\bf Z}\la_N'$, accordingly we have
$P'=P_0'\oplus\cdots \oplus P_{r-1}'$. Likewise we identify
$P$ with $P'_0$. Clearly $P$ and $Q$ are the weight
and root lattices of the root system of types
$C_n$, $B_n$, $2B_n$ and $G_2$ respectively.

Using the same notation as in
section one, the virtual finite dimensional fundamental weights are given by:
\begin{eqnarray*}
\lambda_1&=&\left\{\begin{array}{ll}
     (\lambda_1'+\lambda_{2n-1}')=(\alpha_1'+\cdots
         +\alpha_{2n-1}') &\\
    \ \ =2\alpha_1+\cdots+2\alpha_{n-1}+\alpha_n, &\mbox{for
                 $A_{2n-1}^{(2)}$}\\
     (\lambda_1'+\lambda_2')=(\alpha_1'+\cdots+
          \alpha_{2n}')&\\
     \ \   =2\alpha_1+\cdots+2\alpha_n,
                  &\mbox{for $A_{2n}^{(2)}$}\\
    \la_1'=\a_1'+\cdots+\a_{n-1}'+\frac12(\a_{n}'+\a_{n+1}') &\\
     \ \   =\a_1+\cdots+\a_n, &\mbox{for $D_{n+1}^{(2)}$}\\
    (\la_1'+\la_2'+\la_3')=2(\a_1'+\a_3'+\a_4')+3\a_2'&\\
    \ \  =2\a_1+3\a_2, &\mbox{for $D_4^{(3)}$}
   \end{array}\right. \\
\lambda_n&=&\lambda_n'+\lambda_{n+1}'=\alpha_1'+\cdots
+(n-1)\alpha_{n-1}'+ \frac n2(\alpha_n'+\alpha_{n+1}')\\
&=&\alpha_1+2\alpha_2+\cdots+(n-1)\alpha_n+\frac n2\alpha_n, \ \mbox{for
$D_{n+1}^{(2)}$}.
\end{eqnarray*}
Note that $\Lambda_0$ or $\Lambda_n$ (for $A_{2n}^{(2)}$) is the basic
level one weight.

Let  $ \omega_0 = (-1)^r \omega $. There exists a unique central extension
$$
  1 \longrightarrow
  \langle \omega_0 \rangle
  \longrightarrow
  \hat Q'
  \stackrel{-}{\longrightarrow}
  Q' \longrightarrow 1
$$
of $Q'$ by the cyclic group
  $ \langle \omega_0 \rangle $
with the commutator map $C$ defined below,
\begin{eqnarray*}
  aba^{-1} b^{-1}
&=&
  C (\bar a, \bar b)
  \qquad  \quad
  \mbox{for } a, b \in \hat Q'
\\
&=&
  \prod_{s \in {{\bf Z}}/r{{\bf Z}}}
  (-\omega^s)^{<\bar a|\sigma^s \bar b>}
\end{eqnarray*}

We remark that the commutator map $C$ has the following property:
\begin{eqnarray*}
  &&C (\alpha' + \beta', \gamma') =
   C (\alpha', \gamma') C (\beta', \gamma')
\\
\vspace{-8pt}
  &&C (\alpha', \beta' + \gamma') =
   C (\alpha', \beta') C (\alpha', \gamma')
\\
\vspace{2pt}
  &&C (\alpha', \alpha') = 1
\\
\vspace{-8pt}
  &&C (\sigma \alpha', \sigma \beta') =
   C (\alpha', \beta')
\end{eqnarray*}
for
  $ \alpha' $,
  $ \beta' $,
  $ \gamma' \in Q' $.
Also
  $ C (\alpha', \beta') = C (\beta', \alpha')^{-1} $.

The automorphism
  $ \sigma $
can be lifted to an automorphism
  $ \hat{\sigma} $
of the extension of
  $ \hat Q' $
of $Q'$ such that
\begin{eqnarray}\label{E:char}
  (\hat{\sigma} a)^- &= \sigma \bar a
  \qquad \qquad
  \mbox{for all $ a \in \hat Q'$}
\\
  \hat{\sigma} a &= a
  \qquad \quad
  \mbox{whenever $\sigma \bar a = \bar a $.}
\end{eqnarray}

Let $\a'\to e_{\a'}$ be a section of $\hat Q'$ such that $e_0=1$. Thus
$ e_{\alpha_i'} $,
  $ i=1, \dots , \ell $ are generators of ${\bf C}[\hat Q]$.

We consider the complementary subspace
$$
  Q^{\perp}=
  \left\{
    \alpha \in Q'
    \bigm|
    \langle
       \alpha, \alpha_i
    \rangle = 0, \;
    i \in \Gamma_0
  \right\}=(1-\sigma)Q'
$$
using the diagram automorphism $\sigma$.
By
  $ \hat Q^{\perp} $
we mean the pullback of $Q^{\perp}$ to
  $ \hat Q' $.  Notice that on $Q^{\perp}$ we have $
  C_{Q^{\perp}}(\a, \be)=\omega^{<\sum
  s\sigma^s\a, \be>}$.  Using (\ref{E:char})
we see that there exists a unique
  character $\tau$ : $\hat Q^{\perp} \longrightarrow {\bf C}^{\times}$ such
that
 $$\tau(\omega_0)=\omega_0, \ \ \tau(a\hat\sigma a^{-1})=
 \omega^{-<\sum \sigma^s\overline a|\overline a>/2},  \
\mbox{for } a\in \hat Q'.
 $$

We form the induced $\hat Q'$-module
$$
  {\bf C}\{Q\}=
  {{\bf C}}
  \bigl[
    \hat Q'
  \bigr]
  \otimes_{{\bf C} [\hat Q^{\perp}]}
  {\bf C}_{\tau}
$$
where
${\bf C}_{\tau}$ is the $1$-dimensional
$\hat Q^{\perp}$-module afforded by the character $\tau$. The use of the
notation is justified by the vector space isomorphism:
$V\simeq
 {\bf C}[Q'/Q^{\tau}]$ induced from the natural projection from
$Q'$ to $Q'_{0}=Q$.

On ${\bf C}\{Q\}$ the action of
  $ \hat Q' $
and
  $ Q=
    \{
       \alpha \in Q' \mid
       \sigma (\alpha) = \alpha \} $
are as follows:
\begin{eqnarray*}
  a \cdot b \otimes t
&=&
  a b \otimes t
  \qquad \qquad
  a, b \in \hat Q' , \; t \in T
\\
  \alpha \cdot b \otimes t
&=& <\alpha |\bar b>
  b \otimes t
  \qquad \qquad
  \alpha \in Q .
\end{eqnarray*}

We define
  $ z^{\alpha(0)}
    (\alpha \in Q') $
as an operator on $V$ via
$$
  z^{\alpha (0)}
  \cdot
  b \otimes t =
  z^{\langle \alpha, \bar b \rangle}
  b \otimes t
$$
and the operator
  $ \omega^{\alpha (0)} $
on $U$ by
  $ \omega^{\alpha (0) }
    \cdot
    b \otimes t =
    \omega^{\langle \alpha, \bar b \rangle}
    b \otimes t $.

We have as operators on $V$
\begin{eqnarray*}
  z^{\alpha(0)} a
&=&
  a z^{\alpha(0)+\langle \alpha,\bar a \rangle}
\\
  \omega^{\alpha(0)} a
&=&
  a \omega^{\alpha(0)+\langle \alpha,\bar a \rangle}
\\
  \hat{\sigma} a
&=&
  a \omega^{\Sigma \sigma^r \bar a(0) +
            \langle
                \bar a, \Sigma \sigma^r \bar a
            \rangle / 2 }.
\end{eqnarray*}

For $\a=\sum m_i\a_i'\neq 0$,
we define
$$
e_{\a}=e_{\a_1'}^{m_1}\cdots e_{\a_N'}^{m_N},
$$
where the order is fixed. We then have
\begin{equation}\label{E:3.1}
e_{\a_i'}e_{-\a_i'}=\omega^{\sum<\a_i'|\sigma^s(\a_i')>}, \quad
e^{\a_i'}e^{\a_j'}=(-\omega)^{<\a_i'|\a_j'>}e^{\a_j'}e^{\a_i'}.
\end{equation}
Subsequently we have
\begin{equation}\label{E:3.3}
e_{\a}e_{\be}=\omega^{\sum_{s=0}^{r-1} <\a|\sigma^s (\be)>}e_{\be}e_{\a}.
\end{equation}
for any $\a, \be\in Q'$.

Later we want to use freely $Q$, as fixed point sublattice
of $Q'$, thus we are doing thing not only for $Q'$ but also for $Q$.
For $\la\in P'$ we define the vector space ${\bf C}\{Q\}e_{\la}$ as a
${\bf C}\{Q\}$-module by formally adjoining the element
$e_{\la}$. Then we can define the operator $e^{\lambda}: {\bf C}\{Q\}
->{\bf C}\{Q\}e_{\lambda}$ naturally.

For $\a\in P$ define the operator $\partial_{\alpha}$ on {\bf C}\{Q\}
or ${\bf C}\{Q\}e^{\la}$ by
\begin{equation}\label{E:3.5}
\partial_{\alpha}e^{\be}=(\a|\be)e^{\be}
\end{equation}

Let $U_q(\hat{\bf h})$ be the infinite dimensional Heisenberg algebra
generated by $a_i(k)$ and the central element
$\gamma$ ($k\in {\bf Z}\backslash\{0\}, i=1, \cdots,
n$) subject to the relations
\begin{eqnarray}
a_{\sigma(i)}(k)=a_i\omega^{k}, \label{E:3.6}\\
\ [a_i(k), a_j(l)]&=&\delta_{k+l, 0}\sum_{s=0}^{r-1}
\displaystyle\frac {[<\a_i'|\sigma(\a_j')>/d_i]_i}{k}\omega^{sk}
\displaystyle\frac {\gamma^k-\gamma^{-k}}{q_j-q^{-1}_j}.
\end{eqnarray}
The algebra $U_q(\hat{\bf h})$ acts on the space of symmetric algebra
$Sym(\hat{\bf h}_{\sigma}^-)$ generated by
$a_i(-k)$ ($k\in {\bf N}, i=1, \cdots, n$) with $\gamma =q$
by the action:
\begin{eqnarray*}
&a_i(-k)&\mapsto \mbox{multiplication by}\ a_i(-k),\\
&a_i(k)&\mapsto \sum_{s=0}^{r-1}
\sum_{j}\frac{[<\a_i'|\sigma^s(\a_j')>k/d_i]_i}{k}\omega^{sk}
\frac{q^k-q^{-k}}{q_j-q_j^{-1}}
\frac{d}{d\,a_j(-k)}.
\end{eqnarray*}
For convenience we also identify $Sym(\hat{\bf h}_{\sigma}^-)$ with
$Sym(\hat{\bf h}^-)$ by modulo the relations (\ref{E:3.6}).

We define that
\begin{eqnarray}\nonumber
W_{\epsilon}&=&Sym(\hat{\bf h_{\sigma}}^-)\otimes {\bf C}\{Q\}
, \ \mbox{ for all cases of basic modules}\\
W_1&=&Sym(\hat{\bf h_{\sigma}}^-)\otimes {\bf C}\{Q\}e_{\lambda_1}
, \ \mbox{for $A_{2n-1}^{(2)}$},
\label{E:3.8}\\
W_n&=&Sym(\hat{\bf h}^-)\otimes {\bf C}\{Q\}e_{\la_n}, \
\mbox{ for $D_{n+1}^{(2)}$.}
\nonumber
\end{eqnarray}
where $\epsilon=0$ except for $A_{2n}^{(2)}$, $\epsilon=n$.
Then the algebras $U_q(\hat{\bf h})$ and the space
${\bf C}[P']$ (of course ${\bf C}[P]$)
act on the spaces $W_i$ canonically by extending their
respective actions component-wise. We also define the degree action by
\begin{equation}\label{E:3.9}
d.f\otimes e^{\a}e_{\beta}=(-\sum_{j=1}^km_j-\sum_{j=1}^ln_j-
\frac{<\be|\be>}2 +\sum_{s=0}^{r-1}
\frac{<\alpha_i'|\sigma^s(\alpha_i')>}{2} )f\otimes e^{\a}e_{\beta},
\end{equation}
where
$f\otimes e^{\a}e_{\beta}=
a_{i_1}(-m_{1})\cdots a_{i_k}(-m_k)\otimes e^{\a}e_{\beta}
\in W_i$. For convenience we denote $\la_{\epsilon}=0$.

\begin{prop}\label{P:3.1}
The space $W_i$ ($i=\epsilon, 1, n$) is  the irreducible representation
$V(\Lambda_i)$  of the twisted quantum affine algebra
$U_q(X_N^{(r)})$ under the action:
\begin{eqnarray*}
 \gamma &\mapsto & q, \
K_j \mapsto q^{\frac{p_i}r\partial_{\a_j}},  \  a_j(k)\mapsto a_j(k)
\quad (1\leq j\leq n),\\
x_i^{\pm}(z)&\mapsto& X^{\pm}_i(z)
=\sqrt{p_i}exp(\pm\sum_{k=1}^{\infty}\frac{a_i(-k)}{[k/d_i]_i}q^{\mp k/2}z^k)
 exp(\mp\sum_{k=1}^{\infty}\frac{a_i(k)}{[k/d_i]_i}q^{\mp k/2}z^{-k})\\
& &\hskip 1in \times e^{\pm \a_i}z^{\pm \frac{p_i}r\partial_{\alpha_i}+1},
\end{eqnarray*}
where $p_i=1$ or $r$ according to
$\sigma(i)\neq i$ or $i$, and the degree operator $d$ acts by (\ref{E:3.9}).
The highest weight vectors are respectively:
$$
 |\Lambda_i\rangle
=1\otimes e_{\lambda_i}\otimes 1, \quad i=\epsilon, 1, n.
$$
for appropriate cases.
\end{prop}
{\it Proof}. The detail of the proof is given in \cite{kn:J1} for the case
of $\lambda_{\epsilon}$, and the proof in the cases of $\lambda_1$ for
$A_{2n-1}^{(2)}$ and $\lambda_n$ for $D_{n+1}^{(2)}$ is exactly the same as
that of $\lambda_{\epsilon}$ by the standard shift method of the highest
weight vectors. The irreducibility is easily shown by using the theory of
the
quantum $Z$-algebra operators given by the first author \cite{kn:J2}.
\hskip 1cm $\Box$

\section{Realization of the level one vertex operators}

We first recall the notion of Frenkel-Reshetikhin vertex operators
\cite{kn:FR}.
Let $V$ be a finite dimensional representation of the derived
quantum affine Lie
algebra $U_q'(X_N^{(r)})$ with the associated affinization space $V_z$
(cf. section 2). Let $V(\la)$ and $V(\mu)$ be two integrable highest weight
irreducible representations of $\U$.
The type I (resp. type II) {\it vertex operator} is the  $\U$-intertwining
operator: $V(\la)\longrightarrow V(\mu)\hat{\otimes} V_z$ (resp.
$V(\la)\longrightarrow V_z\hat{\otimes} V(\mu)$). For simplicity we will
compute the
intertwining operators for the derived
subalgebra $U_q'(X_N^{(r)})$, which differ only by a power of $z$ from that
for $\U$. Also we will omit $\hat{\ }$ in the tensor notation.

The existence of vertex operators
is given by the fundamental fact \cite{kn:FR}
(cf. \cite{kn:DJO}) (true for both types, though stated for type I):
\begin{eqnarray}
Hom_{\U}(V(\lambda), V(\mu)\otimes V_z)
&\simeq &\{v\in V|\ wt(v)=\lambda-\mu \ \textstyle{mod} \ \delta
\ \ \mbox{and} \nonumber\\
&& e_i^{\langle \mu, h_i \rangle+1}v=0 \  \mbox {for} \ i=0, \cdots, n\},
\end{eqnarray}
where the isomorphism is defined by sending an element
$\Phi\in Hom_{U_q}(V(\lambda), V(\mu)\otimes V_z)$ to an element $v\in V$ such
that
$$
\Phi|\lambda \rangle=|\mu \rangle \otimes v
+ \mbox{(higher terms in the powers of $z$)}.
$$

For the evaluation module $V_z$, we define the components of
the vertex operator
$\Phi_{\lambda}^{\mu V}:
V(\lambda)\longrightarrow V(\mu)\otimes V_z$ by
\begin{equation} \label{E:4.2}
\Phi_{\lambda}^{\mu V}(z)|u\rangle
=\sum_{j=1}^n\Phi_{\lambda j}^{\mu V}
(z)|u\rangle \otimes v_j+\Phi_{\lambda 0}^{\mu V}(z)|u\rangle\otimes v_0
+\sum_{j=1}^n\Phi_{\lambda \overline{j}}^{\mu V}
(z)|u \rangle \otimes v_{\overline j},
\end{equation}
for $|u \rangle \in V(\lambda)$. The components of type II
vertex operators are defined similarly.

We also consider the intertwining operators of modules of the following form:
$$
\Phi_{\lambda V}^{\mu}(z): V(\lambda)\otimes V_z
\longrightarrow V(\mu)\otimes {\bf C}[z, z^{-1}]
$$
by means of the vertex operators with respect to the dual space $V_z^*$:
\begin{equation}\label{E:4.1}
\Phi_{\lambda V}^{\mu}(z) (|v \rangle
\otimes v_i)=\Phi_{\lambda i}^{\mu V^*}(z)|v \rangle
\end{equation}
for $|v \rangle \in V(\lambda)$ and $i=1, \cdots ,n, 0, \overline{n},
\cdots, \overline{1}$. In other words, the vertex operator
$\Phi_{\lambda i}^{\mu V^*}(z)$ can be defined directly as in the
case of $\Phi_{\lambda i}^{\mu V}(z)$.

>From the above fact we determine that the level one vertex operators exist
only in the following cases:
between
$V(\Lambda_0)$ and $V(\Lambda_{1})$ for $A_{2n-1}^{(2)}$;
between $V(\Lambda_i)$ and itself, $i=0, n$ for $D_{n+1}^{(2)}$;
between $V(\Lambda_n)$ and itself for $A_{2n}^{(2)}$ ,
between $V(\Lambda_0)$ and itself for $D_4^{(3)}$.

We take the following normalization.

For $A_{2n-1}^{(2)}$:
\begin{equation} \label{E:4.3}
\begin{array}{rcl}
\Phi^{\Lambda_{1}V}_{\Lambda_0}(z)|\Lambda_0 \rangle
&=&|\Lambda_{1}\rangle \otimes
v_{\overline{1}}+\mbox{higher terms in $z$},\\
\Phi^{\Lambda_0V}_{\Lambda_{1}}(z)|\Lambda_{1} \rangle
&=&|\Lambda_0 \rangle \otimes
v_{1}+\mbox{higher terms in $z$}\\
\Phi^{\Lambda_{1}V^*}_{\Lambda_0}(z)|\Lambda_0 \rangle
&=&|\Lambda_{1}\rangle \otimes
v_{{1}}^*+\mbox{higher terms in $z$},\\
\Phi^{\Lambda_0V^*}_{\Lambda_{1}}(z)|\Lambda_{1} \rangle
&=&|\Lambda_0 \rangle \otimes
v_{\overline{1}}^*+\mbox{higher terms in $z$};
\end{array}
\end{equation}

For $D_{n+1}^{(2)}$:
\begin{equation}\label{E:4.3'}
\begin{array}{rcl}
\Phi^{\Lambda_{0}V}_{\Lambda_0}(z)|\Lambda_0 \rangle
&=&|\Lambda_{0} \rangle \otimes
v_{\overline 0}+\mbox{higher terms in $z$},\\
\Phi^{\Lambda_{0}V^*}_{\Lambda_0}(z)|\Lambda_0 \rangle
&=&|\Lambda_{0} \rangle \otimes
v_{\overline 0}^*+\mbox{higher terms in $z$}\\
\Phi^{\Lambda_{n}V}_{\Lambda_n}(z)|\Lambda_n \rangle
&=&|\Lambda_{n} \rangle \otimes
v_{0}+\mbox{higher terms in $z$},\\
\Phi^{\Lambda_{n}V^*}_{\Lambda_n}(z)|\Lambda_n \rangle
&=&|\Lambda_{n} \rangle \otimes
v_{0}^*+\mbox{higher terms in $z$}.
\end{array}
\end{equation}

For $A_{2n}^{(2)}$:
\begin{equation}\label{E:4.3''}
\begin{array}{rcl}
\Phi^{\Lambda_{n}V}_{\Lambda_n}(z)|\Lambda_n \rangle
&=&|\Lambda_{n} \rangle \otimes
v_{0}+\mbox{higher terms in $z$},\\
\Phi^{\Lambda_{n}V^*}_{\Lambda_n}(z)|\Lambda_n \rangle
&=&|\Lambda_{n} \rangle \otimes
v_{0}^*+\mbox{higher terms in $z$}.
\end{array}
\end{equation}

For $D_{4}^{(3)}$:
\begin{equation}\label{E:4.3'''}
\begin{array}{rcl}
\Phi^{\Lambda_{0}V}_{\Lambda_0}(z)|\Lambda_0 \rangle
&=&|\Lambda_{0} \rangle \otimes
v_{\overline 0}+\mbox{higher terms in $z$},\\
\Phi^{\Lambda_{0}V^*}_{\Lambda_0}(z)|\Lambda_0 \rangle
&=&|\Lambda_{0} \rangle \otimes
v_{\overline 0}^*+\mbox{higher terms in $z$}.
\end{array}
\end{equation}

Type II vertex operators assume similar normalization.

By the intertwining property it is easy to see the following
determination relations.
\begin{prop}\label{P:4.1}
The vertex operator $\Phi(z)$ of type I with respect to $V_z$ is
determined by its component $\Phi_{\overline 1}(z)$.
With respect to $V_z^*$, the vertex operator $\Phi^*(z)$
of type I is determined by
$\Phi_1^*(z)$.
More explicitly, we have:

For $A_{2n-1}^{(2)}$:
\begin{eqnarray*}
\Phi_i(z)&=&[\Phi_{i+1}(z), f_i]_q
\qquad\mbox{for}\quad i=1, \cdots, n-1,\\
\Phi_{\overline{i+1}}(z)&=&[\Phi_{\overline i}(z), f_i]_q
 \qquad\mbox{for}\quad i=1, \cdots, n-1,\\
\ [f_i, \Phi_i(z)]_q &=&0, \quad [f_i, \Phi_{\overline{i+1}}(z)]_q=0\\
 \ [f_i, \Phi_j(z)]&=&0,  \qquad\mbox{for}\quad j\neq i, i+1, \overline{i+1},
 \overline{i},\\
\Phi_n(z)&=&[\Phi_{\overline n}(z), f_n]_{q^2}, \\
\ [f_n, \Phi_n(z)]_{q^2}&=&0, \quad [f_n, \Phi_j(z)]=0, \quad j\neq n,
\overline n,
\overline{n}.
\end{eqnarray*}

\begin{eqnarray*}
\Phi_{i+1}^*(z)&=&[f_i, \Phi^*_i(z)]_{q^{-1}}
\qquad\mbox{for}\quad i=1, \cdots, n-1,\\
\Phi^*_{\overline{i}}(z)&=&[f_i,
\Phi^*_{\overline{i+1}}(z)]_{q^{-1}}
\qquad\mbox{for}\quad i=1, \cdots, n-1,\\
\ [\Phi_{i+1}^*(z), f_i]_{q^{-1}} &=& 0, \quad
[\Phi_{\overline i}^*(z), f_i]_{q^{-1}}=0,\\
\Phi^*_{\overline n}(z)&=&[f_n, \Phi^*_{n}(z)]_{q^{-2}},\\
\ [\Phi_j^*(z), f_j]&=&0,  \qquad\mbox{for}\quad j\neq i, i+1, \overline{i+1},
 \overline{i},\\
\ [\Phi_n^*(z), f_n]_{q^{-2}}&=&0, \quad [\Phi_j^*(z), f_n]=0,
\quad j\neq n, \overline n,
\overline{n}.
\end{eqnarray*}

For $D_{n+1}^{(2)}$:
\begin{eqnarray*}
\Phi_i(z)&=&[\Phi_{i+1}(z), f_i]_{q^2}
\qquad\mbox{for}\quad i=1, \cdots, n-1,\\
\Phi_{\overline{i+1}}(z)&=&[\Phi_{\overline i}(z), f_i]_{q^2}
 \qquad\mbox{for}\quad i=1, \cdots, n-1,\\
\ [f_i, \Phi_i(z)]_{q^2} &=&0,
\quad [f_i, \Phi_{\overline{i+1}}(z)]_{q^2}=0\\
 \ [f_i, \Phi_j(z)]&=&0,  \qquad\mbox{for}\quad j\neq i, i+1, \overline{i+1},
 \overline{i},\\
\Phi_{n}(z)&=&[\Phi_{0}(z), f_n],\\
\Phi_{0}(z)&=&[2]^{-1}[\Phi_{\overline{n}}(z), f_n]_{q^2}, \\
\ [f_n, \Phi_n(z)]_{q^2}&=&0, \quad [f_n, \Phi_j(z)]=0, \quad j\neq n, 0,
\overline n,\\
\Phi_{\overline 0}(z)&=&[2]^{-1}z[\Phi_1(z), f_0]_{q^2}, \qquad
\Phi_{\overline 1}(z)=z[\Phi_{\overline 0}(z), f_0],\\
\ [\Phi_{\overline 1}(z), f_0]_{q^2}&=& 0, \qquad [\Phi_j(z), f_0]=0,
\quad\mbox{for}\ j\neq 1, \overline 0, \overline 1.
\end{eqnarray*}

\begin{eqnarray*}
\Phi_{i+1}^*(z)&=&[f_i, \Phi^*_i(z)]_{q^{-2}}
\qquad\mbox{for}\quad i=1, \cdots, n-1,\\
\Phi^*_{\overline{i}}(z)&=&[f_i,
\Phi^*_{\overline{i+1}}(z)]_{q^{-2}}
\qquad\mbox{for}\quad i=1, \cdots, n-1,\\
\ [\Phi_{i+1}^*(z), f_i]_{q^{-2}} &=& 0, \quad
[\Phi_{\overline i}^*(z), f_i]_{q^{-2}}=0,\\
\Phi^*_{0}(z)&=&[f_n, \Phi^*_{n}(z)]_{q^{-2}},\quad
\Phi^*_{\overline n}(z)=[2]^{-1}[f_n, \Phi^*_{0}(z)],\\
\ [\Phi_j^*(z), f_i]&=&0,  \qquad\mbox{for}\quad j\neq i, i+1, \overline{i+1},
 \overline{i},\\
\ [f_n, \Phi_{\overline n}^*(z)]_{q^{-2}}&=&0, \quad [\Phi_j^*(z), f_n]=0,
\quad j\neq n, 0, \overline n,
\overline{n}\\
\Phi^*_1(z)&=&[2]^{-1}z[f_0, \Phi^*_{\overline 0}(z)], \quad
\Phi^*_{\overline 0}(z)=z[f_0, \Phi^*_{\overline 1}(z)]_{q^{-2}},\\
\ [\Phi^*_1(z), f_0]_{q^{-2}}&=&0, \quad [\Phi^*_j(z), f_0]=0,
\quad\mbox{for} \ j\neq 1, \overline 0, \overline 1.
\end{eqnarray*}

For $A_{2n}^{(2)}$:
\begin{eqnarray*}
\Phi_i(z)&=&[\Phi_{i+1}(z), f_i]_q
\qquad\mbox{for}\quad i=1, \cdots, n-1,\\
\Phi_{\overline{i+1}}(z)&=&[\Phi_{\overline i}(z), f_i]_q
 \qquad\mbox{for}\quad i=1, \cdots, n-1,\\
\ [f_i, \Phi_i(z)]_q &=&0, \quad [f_i, \Phi_{\overline{i+1}}(z)]_q=0\\
 \ [f_i, \Phi_j(z)]&=&0,  \qquad\mbox{for}\quad j\neq i, i+1, \overline{i+1},
 \overline{i},\\
\Phi_n(z)&=&[\Phi_{0}(z), f_n], \\
\Phi_{0}(z)&=&[2]_n^{-1}[\Phi_{\overline{n}}(z), f_n]_q, \\
\ [f_n, \Phi_n(z)]_{q}&=&0, \quad [f_n, \Phi_j(z)]=0, \quad j\neq n,
\overline n,
\overline{n},\\
\Phi_{\overline 1}(z)&=&z[\Phi_1(z), f_0]_{q^2},\\
\ [f_0, \Phi_{\overline 1}(z)]_{q^2}&=&0 , \quad [f_0, \Phi_j(z)]=0 ,
\quad j\neq 1, \overline 1.
\end{eqnarray*}

\begin{eqnarray*}
\Phi_{i+1}^*(z)&=&[f_i, \Phi^*_i(z)]_{q^{-1}}
\qquad\mbox{for}\quad i=1, \cdots, n-1,\\
\Phi^*_{\overline{i}}(z)&=&[f_i,
\Phi^*_{\overline{i+1}}(z)]_{q^{-1}}
\qquad\mbox{for}\quad i=1, \cdots, n-1,\\
\ [\Phi_{i+1}^*(z), f_i]_{q^{-1}} &=& 0, \quad
[\Phi_{\overline i}^*(z), f_i]_{q^{-1}}=0,\\
\Phi^*_{\overline n}(z)&=&[2]_n^{-1}[f_n, \Phi^*_{0}(z)],\quad
\Phi^*_0(z)=[f_n, \Phi_n^*(z)]_{q^{-1}},\\
\ [\Phi_j^*(z), f_i]&=&0,  \qquad\mbox{for}\quad j\neq i, i+1, \overline{i+1},
 \overline{i},\\
\ [\Phi_n^*(z), f_n]_{q^{-1}}&=&0, \quad [\Phi_j^*(z), f_n]=0,
\quad j\neq n, \overline n,\\
\Phi^*_1(z)&=&z[f-0, \Phi_{\overline 1}^*(z)]_{q^{-2}},\\
\ [\Phi_1^*(z), f_0]_{q^{-2}}&=&0, \quad [\Phi_j^*(z), f_0]_{q^{-2}}=0,
\qquad\mbox{for }\quad j\neq 1, \overline 1.
\overline{n}.
\end{eqnarray*}

For $D_{4}^{(3)}$:
\begin{eqnarray*}
\Phi_1(z)&=&[\Phi_2(z), f_1]_q, \ [f_1, \Phi_1(z)]_q=0,
\ [f_1, \Phi_3(z)]_{q^2}=0, \\
\Phi_3(z)&=&[\Phi_0(z), f_1], \ \Phi_{\overline 0}(z)+[2]\Phi_0(z)=
[\Phi_{\overline 3}(z), f_1]_{q^2},\\
\Phi_{\overline 2}(z)&=&[\Phi_{\overline 1}(z), f_1]_q, \
[f_1, \Phi_{\overline 2}(z)]=0, \ [f_1, \Phi_{\overline 0}(z)]=0, \\
\Phi_2(z)&=&[\Phi_3(z), f_2]_{q^3}, \ \Phi_{\overline 3}(z)=
[\Phi_{\overline 2}(z), f_2]_{q^3},\\
\ [f_2, \Phi_2(z)]_{q^3}&=&0, \ [f_2, \Phi_{\overline 3}(z)]_{q^3}=0,\\
\ [f_2, \Phi_j(z)]&=&0, \quad \mbox{for } j=1, 0, {\overline 1},
{\overline 0}.\\
\ \Phi_{\overline 3}(z)&=& z[\Phi_2(z), f_0]_q, \
\Phi_0(z)+[2]\Phi_{\overline 0}(z)=z[\Phi_1(z), f_0]_{q^2},\\
\Phi_{\overline 2}(z)&=&z[\Phi_3(z), f_0]_q, \ \Phi_{\overline 1}(z)
=z[\Phi_{\overline 0}(z), f_0],\\
\ [f_0, \Phi_j(z)]_q&=&0, \quad\mbox{for }j={\overline 2}, {\overline 3},\\
\ [f_0, \Phi_{\overline 1}(z)]_{q^2}&=&0, \ [f_0, \Phi_0(z)]=0.
\end{eqnarray*}

\begin{eqnarray*}
\Phi_2^*(z)&=&[f_1, \Phi_1^*(z)]_{q^{-1}}, \
[\Phi_{\overline 1}^*(z), f_1]_{q^{-1}}=0,
\ [\Phi_{\overline 3}^*(z), f_1]_{q^{-2}}=0, \\
\Phi_{\overline 3}^*(z)&=&[f_1, \Phi_{\overline 0}^*(z)],
\ [2]\Phi_{\overline 3}^*(z)=
[f_1, \Phi_{0}^*(z)],\\
\Phi_{\overline 1}^*(z)&=&[f_1, \Phi_{\overline 2}^*(z)]_{q^{-1}}, \
[\Phi_{2}^*(z), f_1]_{q^{-1}}=0, \\
\Phi_3^*(z)&=&[f_2, \Phi_2^*(z)]_{q^{-3}}, \ \Phi_{\overline 2}^*(z)=
[f_2, \Phi_{\overline 3}^*(z)]_{q^{-3}},\\
\ [\Phi_{\overline 2}^*(z), f_2]_{q^{-3}}&=&0, \
[\Phi_{ 3}^*(z), f_2]_{q^{-3}}=0,\\
\ [\Phi_j^*(z), f_2]&=&0, \quad \mbox{for } j=1, 0, {\overline 1},
{\overline 0}.\\
\ \Phi_{3}^*(z)&=& z[f_0, \Phi_{\overline 2}(z)]_{q^{-1}}, \
[2]\Phi_{1}^*(z)=z[f_0, \Phi_{\overline 0}^*(z)],\\
\Phi_{2}^*(z)&=&z[f_0, \Phi_{\overline 3}^*(z)]_{q^{-1}},
\ \Phi_{\overline 0}^*(z)
=z[f_0, \Phi_{\overline 1}^*(z)]_{q^{-2}},\\
\ [\Phi_j^*(z), f_0]_{q^{-1}}&=&0, \quad\mbox{for }j={ 2}, {3},\\
\Phi_1^*(z)&=&z[f_0, \Phi_0^*(z)], \ [\Phi_{1}^*(z), f_0]_{q^{-2}}=0.
\end{eqnarray*}

\end{prop}

\begin{prop}\label{P:4.2}
Let $\tilde {\Phi}(z)$ be a type II vertex operator with respect to $V_z$:
$V(\lambda)\longrightarrow V_z\otimes V(\mu)$.
Then $\tilde {\Phi}(z)$ is determined
by the component $\tilde {\Phi}_1(z)$,
and with respect to $V_z^*$ the vertex operator $\Phi^*(z)$
is determined by its component
$\Phi_{\overline{1}}^*(z)$:
 More precisely, with respect to $V_z$,
we have:

For $A_{2n-1}^{(1)}$:
\begin{eqnarray*}
\Phi_{i+1}(z)&=&[\Phi_i(z), e_i]_q \quad i=1, \cdots, n-1,\\
\Phi_{\overline i}(z) &=&[\Phi_{\overline{i+1}}(z), e_i]_q
\quad i=1, \cdots, n-1,\\
\ [e_i, \Phi_{i+1}(z)]_q&=&0, \quad [e_i, \Phi_{\overline{i}}]_q=0,\\
\ [e_i, \Phi_j] &=&0 \qquad \mbox{for}\ j\neq i, i+1, \overline{i+1},
\overline{i},\\
\Phi_{\overline n}(z) &=&[\Phi_{n}(z), e_n]_{q^2},\\
\ [e_n, \Phi_{\overline{n}}(z)]_{q^2} &=& 0, \quad [e_n, \Phi_j(z)]=0 \quad
\mbox{for}\ j\neq n, \overline{n},
\end{eqnarray*}

\begin{eqnarray*}
\Phi_i^*(z)&=&q^2[e_i, \Phi_{i+1}^*(z)]_{q^{-1}}
\quad \mbox{for} \ i=1, \cdots, n-1,\\
\Phi_{\overline{i+1}}^*(z) &=&q^2[e_i,
\Phi_{\overline i}^*(z)]_{q^{-1}} \quad \mbox{for} \ i=1, \cdots, n-1,\\
\ [\Phi_i^*(z), e_i]_{q^{-1}}&=& 0, \quad [\Phi_{\overline{i+1}}^*(z),
e_i]_{q^{-1}}=0, \\
\ [\Phi_j^*(z), e_i]&=& 0 \qquad \mbox{for}\ j\neq i, i+1, \overline{i+1},
\overline{i},\\
\Phi_n^*(z) &=&q^{4}[e_n, \Phi^*_{\overline n}(z)]_{q^{-2}}, \\
\ [\Phi_{n}^*(z), e_n]_{q^{-2}}&=& 0,\quad
[\Phi_j^*(z), e_n]=0 \qquad \mbox{for}\ j\neq n, 0, \overline{n}.
\hskip 1cm \Box
\end{eqnarray*}

For $D_{n+1}^{(2)}$:
\begin{eqnarray*}
\Phi_{i+1}(z)&=&[\Phi_{i}(z), e_i]_{q^2}
\qquad\mbox{for}\quad i=1, \cdots, n-1,\\
\Phi_{\overline{i}}(z)&=&[\Phi_{\overline {i+1}}(z), e_i]_{q^2}
 \qquad\mbox{for}\quad i=1, \cdots, n-1,\\
\ [e_i, \Phi_{i+1}(z)]_{q^2} &=&0,
\quad [e_i, \Phi_{\overline{i}}(z)]_{q^2}=0\\
 \ [e_i, \Phi_j(z)]&=&0,  \qquad\mbox{for}\quad j\neq i, i+1, \overline{i+1},
 \overline{i},\\
\Phi_{\overline n}(z)&=&[\Phi_{0}(z), e_n], \\
\Phi_{0}(z)&=&[2]^{-1}[\Phi_{{n}}(z), e_n]_{q^2}, \\
\ [e_n, \Phi_{\overline n}(z)]_{q^2}&=&0, \quad [e_n, \Phi_j(z)]=0, \quad
j\neq n, 0, \overline n,
\overline{n},\\
\Phi_{\overline 0}(z)&=&[2]^{-1}z^{-1}[\Phi_{\overline 1}(z), e_0]_{q^2},
\qquad
\Phi_{1}(z)=z^{-1}[\Phi_{\overline 0}(z), e_0],\\
\ [e_0, \Phi_{1}(z)]_{q^2}&=& 0, \qquad [\Phi_j(z), e_0]=0,
\quad\mbox{for}\ j\neq 1, \overline 0, \overline 1.
\end{eqnarray*}

\begin{eqnarray*}
\Phi_{i}^*(z)&=&q^4[e_i, \Phi^*_{i}(z)]_{q^{-2}}
\qquad\mbox{for}\quad i=1, \cdots, n-1,\\
\Phi^*_{\overline{i+1}}(z)&=&q^4[e_i,
\Phi^*_{\overline{i}}(z)]_{q^{-2}}
\qquad\mbox{for}\quad i=1, \cdots, n-1,\\
\ [\Phi_{i}^*(z), e_i]_{q^{-2}} &=& 0, \quad
[\Phi_{\overline{i+1}}^*(z), e_i]_{q^{-2}}=0,\\
\Phi^*_{0}(z)&=&q^2[e_n, \Phi^*_{\overline n}(z)],\quad
\Phi^*_{n}(z)=[2]^{-1}[e_n, \Phi^*_{0}(z)],\\
\ [\Phi_j^*(z), e_j]&=&0,  \qquad\mbox{for}\quad j\neq i, i+1, \overline{i+1},
 \overline{i},\\
\ [\Phi_{n}^*(z), e_n]_{q^{-2}}&=&0, \quad [\Phi_j^*(z), e_n]=0,
\quad j\neq n, 0, \overline n,
\overline{n}\\
\Phi^*_{\overline 1}(z)&=&q^2[2]^{-1}z^{-1}[e_0, \Phi^*_{\overline 0}(z)],
\quad
\Phi^*_{\overline 0}(z)=z^{-1}[e_0, \Phi^*_{\overline 1}(z)]_{q^{-2}},\\
\ [\Phi^*_{\overline 1}(z), e_0]_{q^{-2}}&=&0, \quad [\Phi^*_j(z), e_0]=0,
\quad\mbox{for} \ j\neq 1, \overline 0, \overline 1.
\end{eqnarray*}

For $A_{2n}^{(2)}$:
\begin{eqnarray*}
\Phi_{i+1}(z)&=&[\Phi_{i}(z), e_i]_q
\qquad\mbox{for}\quad i=1, \cdots, n-1,\\
\Phi_{\overline{i}}(z)&=&[\Phi_{\overline{i+1}}(z), e_i]_q
 \qquad\mbox{for}\quad i=1, \cdots, n-1,\\
\ [e_i, \Phi_{i+1}(z)]_q &=&0, \quad [e_i, \Phi_{\overline{i}}(z)]_q=0\\
 \ [e_i, \Phi_j(z)]&=&0,  \qquad\mbox{for}\quad j\neq i, i+1, \overline{i+1},
 \overline{i},\\
\Phi_{\overline n}(z)&=&[\Phi_{0}(z), e_n], \\
\Phi_{0}(z)&=&[2]_n^{-1}[\Phi_{{n}}(z), e_n]_q, \\
\ [e_n, \Phi_{\overline n}(z)]_{q}&=&0, \quad [e_n, \Phi_j(z)]=0, \quad
j\neq n, 0,
\overline n,
\overline{n},\\
\Phi_{1}(z)&=&z^{-1}[\Phi_{\overline 1}(z), e_0]_{q^2},\\
\ [e_0, \Phi_{1}(z)]_{q^2}&=&0 , \quad [e_0, \Phi_j(z)]=0 ,
\quad j\neq 1, \overline 1.
\end{eqnarray*}

\begin{eqnarray*}
\Phi_{i}^*(z)&=&q^2[e_i, \Phi^*_{i+1}(z)]_{q^{-1}}
\qquad\mbox{for}\quad i=1, \cdots, n-1,\\
\Phi^*_{\overline{i+1}}(z)&=&q^2[e_i,
\Phi^*_{\overline{i}}(z)]_{q^{-1}}
\qquad\mbox{for}\quad i=1, \cdots, n-1,\\
\ [\Phi_{i}^*(z), e_i]_{q^{-1}} &=& 0, \quad
[\Phi_{\overline {i+1}}^*(z), e_i]_{q^{-1}}=0,\\
\Phi^*_{n}(z)&=&q[2]_n^{-1}[e_n, \Phi^*_{0}(z)],\quad
\Phi^*_0(z)=q[e_n, \Phi_{\overline n}^*(z)]_{q^{-1}},\\
\ [\Phi_j^*(z), e_j]&=&0,  \qquad\mbox{for}\quad j\neq i, i+1, \overline{i+1},
 \overline{i},\\
\ [\Phi_n^*(z), e_n]_{q^{-1}}&=&0, \quad [\Phi_j^*(z), e_n]=0,
\quad j\neq n, \overline n,
\overline{n},\\
\Phi^*_{\overline 1}(z)&=&q^4[e_0, \Phi_1^*(z)]_{q^{-2}},\\
\ [\Phi_{\overline 1}^*(z), e_0]_{q^{-2}}&=&0,
\quad [\Phi_j^*(z), e_0]=0, \quad\mbox{for}\ j\neq 1, \overline 1.
\end{eqnarray*}

For $D_{4}^{(3)}$:
\begin{eqnarray*}
\Phi_{2}(z)&=&[\Phi_{1}(z), e_1]_q , \
[2]\Phi_0(z)+\Phi_{\overline 0}(z)=[\Phi_{3}(z), e_1]_{q^2},\\
\Phi_{\overline 3}(z)&=&[\Phi_0(z), e_1], \
\Phi_{\overline 1}(z)=[\Phi_{\overline 2}(z), e_1]_q, \\
\ [e_1, \Phi_{2}(z)]_q&=&0, \ [e_1, \Phi_{\overline 3}(z)]_{q^2}=0,\\
\ [e_1, \Phi_{\overline 1}(z)]_q&=&0, \ [e_1, \Phi_{\overline 0}(z)]=0\\
\Phi_3(z)&=&[\Phi_2(z), e_2]_{q^3}, \ \Phi_{\overline 2}(z)=
[\Phi_{\overline 3}(z), e_2]_{q^3},\\
\ [e_2, \Phi_j(z)]_{q^3}&=&0, \quad \mbox{for } j=3, {\overline 2},\\
\ [e_2, \Phi_j(z)]&=&0, \quad \mbox{for } j=1, 0 , {\overline 1},
{\overline 0},\\
\Phi_1(z)&=& z^{-1}[\Phi_{\overline 0}(z), e_0], \
\Phi_2(z)=z^{-1}[\Phi_{\overline 3}(z), e_0],\\
\ \Phi_{3}(z)&=& z^{-1}[\Phi_{\overline 2}(z), e_0]_q, \
\Phi_0(z)+[2]\Phi_{\overline 0}(z)=z^{-1}[\Phi_{\overline 1}(z), e_0]_{q^2},\\
\ [e_0, \Phi_j(z)]_q&=&0, \quad\mbox{for }j=2, 3,\\
\ [e_0, \Phi_{1}(z)]_{q^2}&=&0, \ [f_0, \Phi_0(z)]=0.
\end{eqnarray*}

\begin{eqnarray*}
\Phi_{1}^*(z)&=&q^2[e_1, \Phi^*_{2}(z)]_{q^{-1}} , \
\Phi^*_{\overline{2}}(z)=q^2[e_1,
\Phi^*_{\overline{1}}(z)]_{q^{-1}},\\
\Phi_{0}^*(z)&=&q^2[e_1, \Phi^*_{\overline 3}(z)]_{q^{-2}} , \
\Phi^*_{3}(z)=q^2[e_1,
\Phi^*_{0}(z)],\\
\Phi_3^*(z)&=&q^2[e_1, \Phi_{\overline 0}^*(z)], \
\ [\Phi_{1}^*(z), e_1]_{q^{-1}}= 0, \\
\ [\Phi_{3}^*(z), e_1]_{q^{-2}}&=&0, \
[\Phi_{\overline 2}^*(z), e_1]_{q^{-1}}=0,\\
\Phi^*_{2}(z)&=&q^{-6}[e_2, \Phi^*_{3}(z)]_{q^{-3}}, \
\Phi^*_{\overline 3}(z)=q^6[e_2, \Phi_{\overline 2}^*(z)]_{q^{-3}},\\
\ [\Phi_j^*(z), e_2]_{q^{-3}}&=&0,  \qquad\mbox{for}\quad j=2, {\overline 3}\\
\ [\Phi_j^*(z), e_2]&=&0,
\quad\mbox{for } j= 1, 0, \overline 1,
\overline 0.\\
\Phi^*_{\overline 3}(z)&=&z^{-1}q^2[e_0, \Phi_2^*(z)]_{q^{-1}},
\ \Phi^*_{\overline 2}(z)=z^{-1}q^2[e_0, \Phi_3^*(z)]_{q^{-1}},\\
\Phi^*_{\overline 0}(z)&=&z^{-1}q^2[e_0, \Phi_1^*(z)]_{q^{-2}},
\ \Phi^*_{\overline 2}(z)=z^{-1}q^2[e_0,
\Phi_{\overline 0}^*(z)]_{q^{-1} ??},\\
\Phi_{\overline 1}^*(z)&=&z^{-1}q^2[e_0, \Phi_0^*(z)], \
 \ [\Phi_{\overline 1}^*(z), e_0]_{q^{-2}}=0,\\
\ [\Phi_j^*(z), e_0]_{q^{-1}}&=&0, \quad\mbox{for}\ j=
\overline 2, \overline 3.
\end{eqnarray*}

\end{prop}

Next we determine the relations of vertex operators and Drinfeld
generators.
\begin{prop} \label{P:4.3}
(a) Let ${\Phi}(z): V(\lambda)\longrightarrow V(\mu)\otimes V_z$
be a vertex operator of type I,
where  $(\lambda, \mu)=(\Lambda_{0}, \Lambda_1),
(\Lambda_1, \Lambda_0), (\Lambda_{n}, \Lambda_n)$.
Then we have for each $j=1, \cdots, n$ and $k\in {\bf N}$
\begin{eqnarray*}
[{\Phi}_{\overline 1}(z), X_j^+(w)]&=&0,\\
t_j{\Phi}_{\overline 1}(z)t_j^{-1} &=&
q^{\delta_{j1}}{\Phi}_{\overline 1}(z),\\
\ [a_j(k), {\Phi}_{\overline 1}(z)]
      &=&\left\{
      \begin{array}{ll}\delta_{j1}
       \frac{[k]}{k} q^{\frac{4n+3}2 k}(-z)^k{\Phi}_{\overline
1}(z),&\mbox{for}
       \ A_{2n-1}^{(2)}\\
       \delta_{j1}\frac{[l]_1}{l}(-1)^{nl}
       q^{(4n+3)l}z^{2l}{\Phi}_{\overline 1}(z),
       &\mbox{for $k=2l$,}\ D_{n+1}^{(2)}\\
       \delta_{j1}\frac{[k]}{k} q^{\frac{4n+5}2 k}
       (-z)^k{\Phi}_{\overline 1}(z),
       &\mbox{for}\ A_{2n}^{(2)},\\
       \delta_{j1}\frac{[k]}{k} q^{\frac{15}2 k}
       z^k{\Phi}_{\overline 1}(z),
       &\mbox{for}\ D_{4}^{(3)}
       \end{array}\right. \\
       &&\\
\ [a_j(-k), {\Phi}_{\overline 1}(z)]
      &=&\left\{
      \begin{array}{ll}\delta_{j1}
       \frac{[k]}{k} q^{-\frac{4n+1}2 k}(-z)^{-k}{\Phi}_{\overline 1}(z),
       &\mbox{for}\ A_{2n-1}^{(2)}\\
       \delta_{j1}\frac{[l]_1}{l}(-1)^{nl}
       q^{-(2n+1)l}z^{-2l}{\Phi}_{\overline 1}(z),
       &\mbox{for $k=2l$,}\ D_{n+1}^{(2)}\\
       \delta_{j1}\frac{[k]}{k} q^{-\frac{4n+3}2 k}(-z)^{-k}
       {\Phi}_{\overline 1}(z),
       &\mbox{for}\ A_{2n}^{(2)},\\
       \delta_{j1}\frac{[k]}{k} q^{-\frac{13}2 k}
       z^{-k}{\Phi}_{\overline 1}(z),
       &\mbox{for}\ A_{2n}^{(2)}
       \end{array}\right.
\end{eqnarray*}
(b) If ${\Phi}(z)$ is a vertex operator of type I
associated with $V_z^*$, then
\begin{eqnarray*}
[{\Phi}^*_1(z), X^+_j(w)]&=&0,\\
t_j{\Phi}^*_1(z)t_j^{-1}&=&q^{\delta_{j1}}{\Phi}^*_1(z),\\
\ [a_j(k), {\Phi}^*_1(z)]
      &=&\left\{
      \begin{array}{ll}\delta_{j1}
       \frac{[k]}{k} q^{\frac{3k}2}z^k{\Phi}^*_{1}(z),
       &\mbox{for}
       \ A_{2n-1}^{(2)}, A_{2n}^{(2)}, D_4^{(3)}\\
       \delta_{j1}\frac{[l]_1}{l}(-1)^{nl}
       q^{3l}z^{2l}{\Phi}^*_{1}(z),
       &\mbox{for $k=2l$,}\ D_{n+1}^{(2)}
       \end{array}\right. \\
       &&\\
\ [a_j(-k), {\Phi}^*_1(z)]
      &=&\left\{
      \begin{array}{ll}\delta_{j1}
       \frac{[k]}{k} q^{-\frac{k}2}z^{-k}{\Phi}^*_{1}(z),
       &\mbox{for}
       \ A_{2n-1}^{(2)}, A_{2n}^{(2)}, D_4^{(3)}\\
       \delta_{j1}\frac{[l]_1}{l}(-1)^{nl}
       q^{-l}z^{-2l}{\Phi}^*_{1}(z),
       &\mbox{for $k=2l$,}\ D_{n+1}^{(2)}
       \end{array}\right.
\end{eqnarray*}
(c) If ${\Phi}(z)$ is a vertex operator of type II
associated with $V_z$, then
\begin{eqnarray*}
[{\Phi}_1(z), X^-_j(w)]&=&0,\\
t_j{\Phi}_1(z)t_j^{-1}&=&q^{-\delta_{j1}}{\Phi}_1(z),\\
\ [a_j(k), {\Phi}_1(z)]
&=&\left\{
      \begin{array}{ll}-\delta_{j1}
       \frac{[k]}{k} q^{\frac{k}2}z^k{\Phi}_{1}(z),
       &\mbox{for}
       \ A_{2n-1}^{(2)}, A_{2n}^{(2)}, D_4^{(3)}\\
       -\delta_{j1}\frac{[l]_1}{l}(-1)^{nl}
       q^{l}z^{2l}{\Phi}_{1}(z),
       &\mbox{for $k=2l$,}\ D_{n+1}^{(2)}
       \end{array}\right. \\
       &&\\
\ [a_j(-k), {\Phi}_1(z)]&=&\left\{
      \begin{array}{ll}-\delta_{j1}
       \frac{[k]}{k} q^{-\frac{3k}2}z^{-k}{\Phi}_{1}(z),
       &\mbox{for}
       \ A_{2n-1}^{(2)}, A_{2n}^{(2)}, D_4^{(3)}\\
       -\delta_{j1}\frac{[l]_1}{l}(-1)^{nl}
       q^{-3l}z^{2l}{\Phi}_{1}(z),
       &\mbox{for $k=2l$,}\ D_{n+1}^{(2)}
       \end{array}\right.
\end{eqnarray*}
(d) If ${\Phi}^*(z)$ is a vertex operator of type II associated
with $V_z^*$, then
\begin{eqnarray*}
[{\Phi}^*_{\overline 1}(z), X^-_j(w)]&=&0,\\
t_j{\Phi}^*_{\overline 1}(z)t_j^{-1}&=&q^{-\delta_{j1}}
{\Phi}^*_{\overline 1}(z),\\
\ [a_j(k), {\Phi}^*_{\overline 1}(z)]
      &=&\left\{
      \begin{array}{ll}-\delta_{j1}
       \frac{[k]}{k} q^{-\frac{4n-1}2 k}(-z)^k{\Phi}^*_{\overline
1}(z),&\mbox{for}
       \ A_{2n-1}^{(2)}\\
       -\delta_{j1}\frac{[l]_1}{l}(-1)^{nl}
       q^{-(4n-1)l}z^{2l}{\Phi}^*_{\overline 1}(z),
       &\mbox{for $k=2l$,}\ D_{n+1}^{(2)}\\
       -\delta_{j1}\frac{[k]}{k} q^{-\frac{4n+1}2 k}(-z)^k{\Phi}^*_{\overline
1}(z),
       &\mbox{for}\ A_{2n}^{(2)},\\
       -\delta_{j1}\frac{[k]}{k} q^{-\frac{11}2 k}z^k{\Phi}^*_{\overline 1}(z),
       &\mbox{for}\ D_{4}^{(3)}
       \end{array}\right. \\
       &&\\
\ [a_j(-k), {\Phi}^*_{\overline 1}(z)]
      &=&\left\{
      \begin{array}{ll}-\delta_{j1}
       \frac{[k]}{k} q^{\frac{4n-3}2 k}(-z)^{-k}{\Phi}^*_{\overline 1}(z),
       &\mbox{for}\ A_{2n-1}^{(2)}\\
       -\delta_{j1}\frac{[l]_1}{l}(-1)^{nl}
       q^{(2n-3)l}z^{-2l}{\Phi}^*_{\overline 1}(z),
       &\mbox{for $k=2l$,}\ D_{n+1}^{(2)}\\
       -\delta_{j1}\frac{[k]}{k} q^{\frac{4n-1}2 k}(-z)^{-k}{\Phi}^*_{\overline
1}(z),
       &\mbox{for}\ A_{2n}^{(2)},\\
       -\delta_{j1}\frac{[k]}{k} q^{\frac{9}2 k}z^{-k}{\Phi}^*_{\overline
1}(z),
       &\mbox{for}\ D_{4}^{(3)}
       \end{array}\right.
\hskip 1cm \Box
\end{eqnarray*}
\end{prop}

We introduce the auxiliary element $a_{\overline 1}(k) \in
U_q(\hat{\bf h}^-)$ ($k\in {\bf Z}^{\times}$) such that
$$
[a_i(k), a_{\overline{1}}(l)]=\delta_{i1}\delta_{k, -l},
$$
where $k$ and $l$ are even integers in the case of $D_{n+1}^{(2)}$.
We determine that
\begin{eqnarray*}
A_{2n-1}^{(2)}: \qquad
a_{\overline 1}(k)&=&\frac{-k}{[k]^2}
( \sum_{i=1}^{n-1}\frac{q^{(n-i)k}+(-1)^kq^{-(n-i)k}}
{q^{nk}+(-1)^kq^{-nk}}a_i(k)\\
&&\qquad+\frac{q+q^{-1}}{q^{nk}+q^{-nk}}a_n(k) ), \\
D_{n+1}^{(2)}: \qquad
a_{\overline 1}(2k)&=&\frac{-k}{[k]_1^2}
( \sum_{i=1}^{n-1}\frac{q^{2(n-i)k}+q^{-2(n-i)k}}
{q^{2nk}+q^{-2nk}}a_i(2k)\\
&&\qquad+\frac{2}{[2](q^{2nk}+q^{-2nk})}a_n(2k) )\\
A_{2n}^{(2)}: \qquad
a_{\overline 1}(k)&=&\frac{-k}{[k]^2}
( \sum_{i=1}^{n-1}\frac{q^{(2n-2i+1)k/2}+(-1)^kq^{-(2n-2i+1)k/2}}
{q^{(2n+1)k/2}+(-1)^kq^{-(2n+1)k/2}}a_i(k)\\
&&\qquad+\frac{q^{k/2}+(-1)^kq^{-k/2}}{(q^{1/2}+q^{-1/2})
(q^{(2n+1)k/2}+(-1)^kq^{-(2n+1)k/2})}a_n(k) )\\
D_4^{(3)}: \qquad
a_{\overline 1}(k)&=&\frac{-k}{[k][2k]}
( \frac{\prod_{i=0}^{2}(q^{k}+\omega^{ik}q^{-k})}
{q^{3k}+q^{-3k}}a_1(k)\\
&&\qquad+\frac{(q^k+q^{-k})^2}{q^{3k}+q^{-3k}}a_2(k) ).
\end{eqnarray*}
where $a_n(k)=0$ if $k\neq 0\ mod\ r$ for $A_{2n-1}^{(2)}$ and $D_4^{(3)}$
by definition.

We now state our main results on the realization of vertex operators
for the level one $\U$-modules. We only need to
determine one component for each vertex operator based on the previous
propositions.

\begin{theorem} The $\overline 1$-components of the type I vertex operator
${\Phi}(z)_{\lambda}^{\mu V}$ with respect to $V_z$
$: V(\lambda)\longrightarrow V(\mu)\otimes V_z$ can be realized explicitly as
follows:
\begin{eqnarray*}
{\Phi}_{\overline 1}(z)
&=exp(\sum\frac{[k]}kq^{\frac{4n+3}2k}a_{\overline 1}(-k)(-z)^k)
exp(\sum\frac{[k]}kq^{-\frac{4n+1}2k}a_{\overline 1}(k)(-z)^{-k})\\
&\quad \times e^{\la_1}(-q^{2n+1}z)^{\partial_{\la_1}+(\la_1|\la_1-
\overline{\mu})}
b_{\la}^{\mu},\ \mbox{for } A_{2n-1}^{(2)},\\
{\Phi}_{\overline 1}(z)
&=exp(\sum\frac{[k]_1}kq^{(4n+3)k}(-1)^{kn}a_{\overline 1}(-2k)z^{2k})
exp(\sum\frac{[k]_1}kq^{-(2n+1)k}(-1)^{kn}a_{\overline 1}(2k)z^{-2k})\\
&\quad \times e^{\la_1}(-q^{2n+1}z)^{\partial_{\la_1}+(\la_1|\la_1-
\overline{\mu})}
b_{\la}^{\mu},\ \qquad\mbox{for } D_{n+1}^{(2)},\\
{\Phi}_{\overline 1}(z)
&=exp(\sum\frac{[k]}kq^{\frac{4n+5}2k}a_{\overline 1}(-k)(-z)^k)
exp(\sum\frac{[k]}kq^{-\frac{4n+3}2k}a_{\overline 1}(k)(-z)^{-k})\\
&\quad \times e^{\la_1}(-q^{2n+2}z)^{\partial_{\la_1}+(\la_1|\la_1-
\overline{\mu})}
b_{\la}^{\mu},\ \qquad\mbox{for } A_{2n}^{(2)},\\
{\Phi}_{\overline 1}(z)
&=exp(\sum\frac{[k]}kq^{\frac{15}2k}a_{\overline 1}(-k)z^k)
exp(\sum\frac{[k]}kq^{-\frac{13}2k}a_{\overline 1}(k)z^{-k})\\
&\quad \times e^{\la_1}(-q^{7}z)^{\partial_{\la_1}+(\la_1|\la_1-
\overline{\mu})}
b_{\la}^{\mu},\ \qquad\mbox{for } D_{4}^{(3)},
\end{eqnarray*}
where $b_{\la}^{\mu}$ is a constant for each case ($b_{\Lambda_0}^{\mu}$=1)
and
$-$ is the canonical projection: $\hat{P}\longrightarrow P$.

The $1$-component of type I vertex operator $\Phi_{\la}^{\mu V^*}(z)$
associated with $V_z^*$ is given by
\begin{eqnarray*}
{\Phi}_1^*(z)&=&exp(\sum\frac{[k]}kq^{3k/2}a_{\overline 1}(-k)z^k)
exp(\sum\frac{[k]}kq^{-k/2}a_{\overline 1}(k)z^{-k})\\
\ && \quad\times e^{\la_1}(-qz)^{\partial_{\la_1}+(\la_1|\la_1-
\overline{\mu})}
b_{\la}^{\mu}, \qquad\mbox{for } A_{2n-1}^{(2)}, A_{2n}^{(2)}, D_4^{(3)}\\
{\Phi}_1^*(z)&=&exp(\sum\frac{[k]}kq^{3k/2}(-1)^{kn}a_{\overline 1}(-2k)z^k)
exp(\sum\frac{[k]}kq^{-k/2}(-1)^{kn}a_{\overline 1}(2k)z^{-2k})\\
\ && \quad\times e^{\la_1}(qz)^{\partial_{\la_1}+(\la_1|\la_1-
\overline{\mu})}
b_{\la}^{\mu}, \qquad\mbox{for } D_{n+1}^{(2)},
\end{eqnarray*}
where $b_{\la}^{\mu}$ is a constant for each case ($b_{\Lambda_0}^{\mu}$=1).
\end{theorem}
{\it Proof.} The theorem is shown by demonstrating that
the construction satisfies
all the intertwining relations.

The exponential factors of Heisenberg generators are required by the
commutation relations in Proposition \ref{P:4.3} and the
definition of $a_{\overline 1}(k)$. The commutation relation of $\Phi(z)$
or $\Phi(z)^*$ with $t_j$ assert the factor $e^{\la_1}$. The various factors
involving $\partial_{\la_1}$ simply assure that the vertex operators
commute with $X^+_j(z)$. For instance in the case of type I vertex operator
with respect to $V_z$ we have
\begin{eqnarray*}
\ [\Phi_{\overline 1}(z), X_j^+(w)]&=&:\Phi_{\overline 1}(z)X_j^+(w):
\ka(z)^{\delta_{jn}}(q^{2n}z)^{(\la_1|\la_1-\overline{\mu})}
(-1)^{(\la_n|\la_1-\overline{\mu})}w^{(\alpha_j|\alpha_j)/2}\\
&&((q^{2n}z-w)^{\delta_{j1}}e^{\la_1}e^{\alpha_j}(-1)^{\delta_{jn}}-
(w-q^{2n}z)^{\delta_{j1}}e^{\alpha_j}e^{\la_1})
(q^{2n}z)^{\partial_{\la_1}}
w^{\partial_{\alpha_j}}(-1)^{\partial_{\la_n}}\\
&=&0.
\end{eqnarray*}
where $:\quad:$ is the bosonic normal ordering for the Heisenberg generators.

The appearance of $(\a_1|\la_1-\overline{\mu})$ and $b_{\la}^{\mu}$ are due to
our
normalization, which is illustrated in the case:
$V(\Lambda_1) \longrightarrow V(\Lambda_0)\otimes V_z$ for $A_{2n-1}^{(2)}$.

\begin{eqnarray*}
{\Phi}_{\Lambda_1\overline 0}^{\Lambda_0 V}(0)|\Lambda_1\rangle
&=&\Phi_{\overline 1}(0)X^-_1(0)\cdots X_n^-(0)\cdots X_1^-(0).1
\otimes e^{\la_1}
\otimes 1 +\cdots\\
&=& b_{\la_1}^{\la_0}e^{\la_1}(-q^{2n+1})^
{\partial_{\la_1}+(\la_1|\la_1-\la_0)}.1\otimes
e^{-\alpha_1}\cdots e^{-\alpha_n}\cdots e^{-\a_1}e^{\lambda_1}\otimes 1 \\
&=&|\Lambda_0 \rangle.
\end{eqnarray*}
\hfill $\Box$
\medskip

By the same argument,
we get the realization for the vertex operators of type II.

\begin{theorem} The $1$-components of the type II vertex operator
${\Phi}_{\lambda}^{V\mu}(z)$ with respect to $V_z$
$: V(\lambda)\longrightarrow V_z\otimes V(\mu)$ can be realized explicitly as
follows:
\begin{eqnarray*}
{\Phi}_{1}(z)
&=&exp(\sum-\frac{[k]}kq^{\frac k2}a_{\overline 1}(-k)z^k)
exp(-\sum\frac{[k]}kq^{-\frac{-3k}2}a_{\overline 1}(k)z^{-k})\\
\ &&\quad \times e^{-\la_1}(-qz)^{-\partial_{\la_1}+(\la_1|\la_1+
\overline{\mu})}
b_{\la}^{\mu}, \qquad\mbox{for } A_{2n-1}^{(2)}, A_{2n}^{(2)}, D_4^{(3)}\\
&=&exp(\sum-\frac{[k]_1}kq^{k}(-1)^{kn}a_{\overline 1}(-2k)z^{2k})
exp(-\sum\frac{[k]_1}kq^{-3k}a_{\overline 1}(2k)z^{-2k})\\
\ &&\quad \times e^{-\la_1}(-qz)^{-\partial_{\la_1}+(\la_1|\la_1+
\overline{\mu})}
b_{\la}^{\mu}, \qquad\mbox{for } D_{n+1}^{(2)},\\
\end{eqnarray*}
where $b_{\la}^{\mu}$ is as above.
The $\overline{1}$-component of type I vertex operator
$\Phi_{\overline{1}}^*(z)$
associated with $V_z^*$ is given by
\begin{eqnarray*}
{\Phi}_{\overline{1}}^*(z)
&=exp(-\sum\frac{[k]}kq^{{\frac{-4n+1}2}k}a_{\overline 1}(-k)(-z)^k)
exp(-\sum\frac{[k]}kq^{\frac{4n-3}2k}a_{\overline 1}(k)(-z)^{-k})\\
&\quad \times e^{-\la_1}(q^{-2n+2}z)^{-\partial_{\la_1}+(\la_1|\la_1+
\overline{\mu})}
b_{\la}^{\mu}, \qquad\mbox{for } A_{2n-1}^{(2)}\\
{\Phi}_{\overline{1}}^*(z)
&=exp(-\sum\frac{[k_1]}k(-1)^{kn}q^{(-4n+1)k}a_{\overline 1}(-2k)z^{2k})
exp(-\sum\frac{[k]}k(-1)^{kn}q^{(4n-5)k}a_{\overline 1}(k)z^{-2k})\\
&\quad \times e^{-\la_1}(q^{-2n+2}z)^{-\partial_{\la_1}+(\la_1|\la_1+
\overline{\mu})}
b_{\la}^{\mu}, \qquad\mbox{for } D_{n+1}^{(2)}\\
{\Phi}_{\overline{1}}^*(z)
&=exp(-\sum\frac{[k]}kq^{{\frac{-4n-1}2}k}a_{\overline 1}(-k)(-z)^k)
exp(-\sum\frac{[k]}kq^{\frac{4n-1}2k}a_{\overline 1}(k)(-z)^{-k})\\
&\quad \times e^{-\la_1}(q^{-2n}z)^{-\partial_{\la_1}+(\la_1|\la_1+
\overline{\mu})}
b_{\la}^{\mu}, \qquad\mbox{for } A_{2n}^{(2)}\\
{\Phi}_{\overline{1}}^*(z)
&=exp(-\sum\frac{[k]}kq^{{\frac{-11}2}k}a_{\overline 1}(-k)z^k)
exp(-\sum\frac{[k]}kq^{\frac{9}2k}a_{\overline 1}(k)z^{-k})\\
&\quad \times e^{-\la_1}(q^{-5}z)^{-\partial_{\la_1}+(\la_1|\la_1+
\overline{\mu})}
b_{\la}^{\mu}, \qquad\mbox{for } D_{4}^{(3)}.
\end{eqnarray*}
\end{theorem}
\hfill $\Box$

\end{document}